\documentclass[preprint]{elsarticle}
\usepackage[utf8]{inputenc}

\usepackage{amsmath,amsbsy,amsthm,amssymb}
\usepackage{multirow,multicol,tabularx,booktabs}
\usepackage{graphicx,placeins,color,url}
\usepackage{bbm,bm}
\usepackage[ruled]{algorithm2e}
\usepackage[shortlabels]{enumitem}

\usepackage{mathrsfs}

\usepackage[top=2.5cm,bottom=2.5cm,right=2.5cm,left=2.5cm]{geometry}
\usepackage[hidelinks]{hyperref}

\usepackage{setspace} 
\newtheorem{conj}{Conjecture}
\newtheorem{thm}[conj]{\bf Theorem}

\newtheorem{cor}[conj]{\bf Corollary}
\newtheorem{prop}[conj]{\bf Proposition}
\newtheorem{lemma}[conj]{\bf Lemma}
\newtheorem{rem}[conj]{\bf Remark}

\newtheorem{assumpt}{\bf Assumption}

\def\bar{\overline}
\def\sgn{{\rm\ sign}}

\def\implies{\Longrightarrow}

\def\to{\rightarrow}

\def\inprobto{{\buildrel \mathbb P \over \longrightarrow \,}}
\def\indistrto{\buildrel {D} \over \longrightarrow}

\def\simiid{\buildrel {\text{i.i.d.}} \over \sim}

\def\Ac{\mbox{$\mathcal A$}}

\def\Nc{\mbox{$\mathcal N$}}

\def\Sc{\mathcal S}
\def\Tc{\mbox{$\mathcal T$}}

\def\Xc{\mbox{$\mathcal X$}}

\def\Zc{{\mathcal Z}}

\def\FF{{\mathbb F}}
\def\HH{{\mathbb H}}

\def\GG{{\mathbb G}}

\def\Rb{\mbox{$\mathbb R$}}

\def\Zb{\mbox{$\mathbb Z$}}

\newcommand{\Ifrak}[2]{{\mathfrak{I}_{#1,#2}}}
\def\Inn{{\mathfrak{I}_{n,n}}}
\def\Ikn{{\mathfrak{I}_{k,n}}}
\def\Iknpr{{\mathfrak{I}_{k,n'}}}
\def\Iknup{{\mathfrak{I}_{k,n}^{\uparrow}}}

\def\EE{ {\rm I} \kern-.15em {\rm E} }
\def\PP{ {\rm I} \kern-.15em {\rm P} }

\def\t{ {\bf t}}
\def\u{ {\bf u}}
\def\v{ {\bf v}}
\def\x{ {\bf x}}
\def\y{ {\bf y}}
\def\z{ {\bf z}}

\def\W{ {\bf W}}
\def\X{ {\bf X}}
\def\Y{ {\bf Y}}
\def\Z{ {\bf Z}}
\def\psibm{ {\bm \psi}}

\def\mds{\medskip}
\def \1{\mathbbm{1} }

\begin{document}

\begin{frontmatter}

    \title{Estimation of a regular conditional functional by \\ conditional U-statistics regression}
    
    \author[crest]{Alexis Derumigny}
    \ead{alexis.derumigny@ensae.fr}
    \address[crest]{CREST-ENSAE, 5, avenue Henry Le Chatelier,
    91764 Palaiseau Cedex, France.}
    
    \date{\today}
    
    \begin{abstract}
        U-statistics constitute a large class of estimators, generalizing the empirical mean of a random variable $\X$ to sums over every $k$-tuple of distinct observations of $\X$. They may be used to estimate a regular functional $\theta(\PP_{\X})$ of the law of $\X$. When a vector of covariates $\Z$ is available, a conditional U-statistic may describe the effect of $\z$ on the conditional law of $\X$ given $\Z=\z$, by estimating a regular conditional functional $\theta(\PP_{\X|\Z=\cdot})$.
        We prove concentration inequalities for conditional U-statistics.
        Assuming a parametric model of the conditional functional of interest, we propose a regression-type estimator based on conditional U-statistics.
        Its theoretical properties are derived, first in a non-asymptotic framework and then in two different asymptotic regimes.
        Some examples are given to illustrate our methods. \\
    \end{abstract}


\end{frontmatter}



%
%
%

\section{Introduction}

Let $\X$ be a random element with values in a measurable space $(\Xc, \Ac)$, and denote by $\PP_{\X}$ its law. A natural framework is $\Xc = \Rb^{p_X}$, for a fixed dimension $p_{\X} > 0$. Often, we are interested in estimating a regular functional $\theta(\PP_{\X})$ of the law of $\X$, of the form
$$\theta(\PP_{\X})
    = \EE \big[ g(\X_1, \dots, \X_k) \big]
    = \int g(\x_1, \dots, \x_k) d\PP_{\X}(\x_1) \cdots d\PP_{\X}(\x_k),
$$
for a fixed $k > 0$, a function $g : \Xc^k \to \Rb$ and $\X_1, \dots \X_k \simiid \PP_{\X}$.
Following Hoeffding~\cite{hoeffding1948class}, a natural estimator of $\theta(\PP_{\X})$ is the U-statistics $\hat \theta(\PP_{\X})$, defined by
\begin{align*}
    \hat \theta(\PP_{\X})
    := |\Ikn|^{-1} \sum_{\sigma \in \Ikn}
    g\big(\X_{\sigma(1)}, \dots, \X_{\sigma(k)} \big),
\end{align*}
where $\Ikn$ is the set of injective functions from $\{1, \dots, k\}$ to $\{1, \dots, n\}$.
For an introduction to the theory of U-statistics, we refer to Koroljuk and Borovskich~\cite{koroljuk1994theory} and Serfling~\cite[Chapter 5]{serfling1980approximation}

\mds

In our framework, we assume that we actually observe $(\X, \Z)$ where $\Z$ is a $p$-dimensional covariate. We are now interested in regular functionals of the conditional law $\PP_{\X | \Z}$. For each $\z_1, \dots, \z_k \in \Zc$, where $\Zc$ is a compact subset of $\Rb^p$, we can define such a functional $\theta_{\z_1, \dots, \z_k}$ by
\begin{align*}
    \theta_{\z_1, \dots \z_k}( \PP_{\X | \Z = \cdot})
    &:= \theta( \PP_{\X | \Z = \z_1}, \dots, \PP_{\X | \Z = \z_k} ) \\
    &= \EE_{\bigotimes_{i=1}^k \PP_{\X | \Z = \z_i}}
    \big[ g(\X_1, \dots, \X_k) \big]
    = \EE \big[ g(\X_1, \dots, \X_k) \big| \Z_i=\z_i, \forall i=1, \dots, k \big] \\
    &= \int g(\x_1, \dots, \x_k) d\PP_{\X | \Z = \z_1}(\x_1) \cdots d\PP_{\X | \Z = \z_k}(\x_k).
\end{align*}
This can be seen as a generalization of $\theta(\PP_\X)$ to the conditional case. Indeed, when $\X$ and $\Z$ are independent,
the new functional $\theta_{\z_1, \dots, \z_k}( \PP_{\X | \Z = \cdot})$ is equal to the unconditional functional $\theta(\PP_\X)$.
For convenience, we will use the notation
$\theta(\z_1, \dots, \z_k)
:= \theta_{\z_1, \dots \z_k}( \PP_{\X | \Z = \cdot}),$
treating the law of $(\X, \Z)$ as fixed (but unknown).
Stute~\cite{stute1991conditional} defined a kernel-based estimator $\hat \theta(\z_1, \dots, \z_k)$ of the conditional functional $\theta(\z_1, \dots, \z_k)$ by
\begin{equation}
    \hat \theta(\z_1, \dots, \z_k)
    := \dfrac{\sum_{\sigma \in \Ikn} K_h \big(\Z_{\sigma(1)}-\z_1 \big) \cdots
    K_h \big(\Z_{\sigma(k)}-\z_k \big) \, g\big(\X_{\sigma(1)}, \dots, \X_{\sigma(k)} \big) }
    {\sum_{\sigma \in \Ikn} K_h \big(\Z_{\sigma(1)}-\z_1 \big) \cdots K_h \big(\Z_{\sigma(k)}-\z_k \big) },
    \label{condUstat:def:estimator_hat_theta}
\end{equation}
where $h > 0$ is the bandwidth, $K(\cdot)$ a kernel on $\Rb^p$, $K_h(\cdot) := h^{-p} K(\cdot / h)$, and $(\X_i, \Z_i) \simiid \PP_{\X, \Z}$.
Stute~\cite{stute1991conditional} proved the asymptotic normality of $\hat \theta(\z_1, \dots, \z_k)$ and its weak and strong consistency. Dony and Mason~\cite{dony2008uniform} derived its uniform in bandwidth consistency under VC-type conditions over a class of possible functions $g$.

\mds

Nevertheless, the estimator (\ref{condUstat:def:estimator_hat_theta}) has several weaknesses. First, the interpretation of the whole hypersurface $(\z_1, \dots, \z_k) \mapsto \hat \theta(\z_1, \dots, \z_k)$ can be difficult. Indeed, the latter curve is of dimension $1+p \times k$, and it is rather challenging to visualize it even for small values of $p$ and $k$.
Second, for each new $k$-uplet $(\z_1, \dots, \z_k)$, the computation of $\hat \theta(\z_1, \dots, \z_k)$ has a cost of $O(n^k)$.
Then, if we want to estimate $\hat \theta(\z_1^{(i)}, \dots, \z_k^{(i)})$ for every $i=1, \dots, N$,
where
$\big(\z_1^{(1)}, \dots, \z_k^{(1)}, \dots,
\z_1^{(N)}, \dots, \z_k^{(N)} \big)
\in \Zc^{k \times N}$, then the total cost is
$O(N n^k)$.
Third, it is well-known that kernel estimators are not very smooth, in the sense that they usually present many spurious local minima and maxima, and this can be a problem in some applications. Therefore, we may want to build estimators which are more regular with respect to the conditioning variables $\z_1, \dots \z_k$, and have a simple functional form.

\mds

Another idea is to decompose the function $(\z_1, \dots, \z_k) \mapsto \theta(\z_1, \dots, \z_k)$ on a basis $(\psi_i)_{i \geq 0}$, generalizing the work of Derumigny and Fermanian \cite{derumigny2018kendall}.
This may not be always easy if the range of the function $\theta(\cdot, \cdots, \cdot)$ is a strict subset of $\Rb$. In that case, it is always possible to use a ``link function''~$\Lambda$, that would be strictly increasing and continuously differentiable and such that the range $\Lambda \circ \theta(\cdot, \cdots, \cdot)$ is exactly $\Rb$.
Whatever the choice of $\Lambda$ (including the identity function), we can decompose the latter function on any basis $(\psi_i)_{i \geq 0}$.
If only a finite number $r > 0$ of elements of this basis are necessary to represent the whole function $\Lambda \circ \theta(\cdot, \cdots, \cdot)$ over $\Zc^k$, then we have the following parametric model:
\begin{equation}
    \forall (\z_1, \dots, \z_k) \in \Zc^k, \,
    \Lambda \big( \theta(\z_1, \dots, \z_k) \big) = \psibm (\z_1, \dots, \z_k)^T \beta^*,
    \label{model:lambda_cond_theta_Z}
\end{equation}
where $\beta^* \in \Rb^r$ is the true parameter
and $\psibm(\cdot) := \big(\psi_1(\cdot), \dots, \psi_r(\cdot) \big)^T \in \Rb^r$.
In most applications, finding an appropriate basis $\psibm$ is not easy. This will depend on the choice of the (conditional) functional $\theta$.
Therefore, the most simple solution consists in choosing a concatenation of several well-known basis such as polynomials, exponentials, sinuses and cosinuses, indicator functions, etc... 
They allow to take into account potential non-linearities and even discontinuities of the function $\Lambda \circ \theta(\cdot, \cdots, \cdot)$.
For the sake of inference, a necessary condition is the linear independence of such functions, as seen in the following proposition (whose straightforward proof is omitted).
\begin{prop}
    The parameter $\beta^*$ is identifiable in Model (\ref{model:lambda_cond_theta_Z}) if and only if the functions
    $(\psi_1(\cdot), \dots, \psi_r (\cdot))$ are linearly independent $\PP_\Z^{\otimes n}$-almost everywhere in the sense that, for all vectors $\t = (t_1, \dots, t_r) \in \Rb^r$,
    $\PP_\Z^{\otimes n}
    \big( \psibm(\Z_1, \dots, \Z_n)^T \t = 0 \big) = 1 \implies \t = 0.$
    \label{condUstat:prop:identifiability_condition}
\end{prop}




With such a choice of a wide and flexible class of functions, it is likely that not all these functions are relevant.
This is know as sparsity, i.e. the number of non-zero coefficients of $\beta^*$, denoted by $|\Sc|$ = $|\beta^*|_0$ is less than $s$, 
for some $s \in \{1, \dots, r\}$. Here, 
$| \cdot |_0$ denotes the number of non-zero components of a vector of $\Rb^{r}$ and 
$\Sc$ is the set of non-zero components of $\beta^*$.
Note that, in this framework, $r$ can be moderately large, for example $30$ or $50$, while the original dimension $p$ is small, for example $p=1$ or~$2$.
This corresponds to the decomposition of a function, defined on a small-dimension domain, in a mildly large basis.

\begin{rem}
	At first sight, in Model (\ref{model:lambda_cond_theta_Z}), there seems to be no noise perturbing the variable of interest. In fact, this can be seen as a simple consequence of our formulation of the model. In the same way, the classical linear model $Y=\X^T \beta^* + \varepsilon$ can be rewritten as $\EE[Y|\X=\x] = \x^T \beta^*$ without any explicit noise. By definition, $\EE[Y|\X=\x]$ is a deterministic function of a given $\x$. In our case, the corresponding fact is: $\Lambda \big( \theta(\z_1, \dots, \z_k) \big)$ is a deterministic function of the variables $(\z_1, \dots, \z_k)$. This means that we cannot write formally a model with noise, such as $\Lambda \big( \theta(\z_1, \dots, \z_k) \big) = \psibm (\z_1, \dots, \z_k)^T \beta^* + \varepsilon$ where $\varepsilon$ is independent of the choice of $(\z_1, \dots, \z_k)$ since the left-hand side of the latter equality is a $(\z_1, \dots, \z_k)$-mesurable quantity, unless $\varepsilon$ is constant almost surely.
\end{rem}

Contrary to more usual models, the explained variable
$\Lambda \big( \theta(\z_1, \dots, \z_k) \big)$, is not observed in Model~(\ref{model:lambda_cond_theta_Z}).
Therefore, a direct estimation of the parameter $\beta^*$ (for example, by the ordinary least squares, or by the Lasso) is unfeasible.
In other words, even if the function $(\z_1, \dots, \z_k) \mapsto \Lambda \big( \theta(\z_1, \dots, \z_k) \big)$ is deterministic (by definition of conditional probabilities), finding the best $\beta$ in Model (\ref{model:lambda_cond_theta_Z}) is far from being a numerical analysis problem since the function to be decomposed is unknown.
Nevertheless, we will replace $\Lambda \big( \theta(\z_1, \dots, \z_k) \big)$ by the nonparametric estimate $\Lambda \big( \hat \theta(\z_1, \dots, \z_k) \big)$, and use it as an approximation of the explained variable.

\mds

More precisely, we fix a finite collection of points $\z'_1, \dots, \z'_{n'} \in \Zc^{n'}$ and a collection $\Iknpr$ of injective functions $\sigma: \{1, \dots, k\} \to \{1, \dots, n'\}$. Note that we are not forced to include \textit{all} the injective functions in $\Iknpr$, reducing its number of elements. This will allow us to decrease the computational cost of the procedure.
For every $\sigma \in \Iknpr,$ we estimate $\hat \theta(\z'_{\sigma(1)}, \dots, \z'_{\sigma(k)})$.
Finally, the estimator $\hat \beta$ is defined as the minimizer of the following $l_1$-penalized criteria
\begin{equation}
    \hat \beta := \arg \min_{\beta \in \Rb^{r}}
    \left[ \frac{(n'-k)!}{n'!}
    \sum_{\sigma \in \Iknpr} \bigg( \Lambda
    \Big(\hat \theta \big(\z'_{\sigma(1)}, \dots, \z'_{\sigma(k)} \big) \Big)
    - \psibm \big(\z'_{\sigma(1)}, \dots, \z'_{\sigma(k)} \big)^T \beta \bigg)^2 + \lambda |\beta|_1 \right],
    \label{condUstat:def:estimator_hat_beta}
\end{equation}
where $\lambda$ is a positive tuning parameter (that may depend on $n$ and $n'$), and $|\cdot|_q$ denotes the $l_q$ norm, for $1\leq q \leq \infty$. This procedure is summed up in the following Algorithm~\ref{condUstat:algo:estimation_beta}.
Note that even if we study the general case with any $\lambda \geq 0$, the corresponding properties of the unpenalized estimator can be derived by choosing the particular case $\lambda = 0$.

\begin{algorithm}[htb]
\label{condUstat:algo:estimation_beta}
\SetAlgoLined
    \vspace{0.1cm}
    \KwIn{A dataset $(X_{i,1}, X_{i,2}, \Z_i)$, $i=1,\dots,n$}
    \KwIn{A finite collection of points $\z'_1, \dots, \z'_{n'} \in \Zc^{n'}$, selected for estimation}
    \KwIn{A collection of $N$ $k$-tuples for prediction $\big(\z_1^{(1)}, \dots, \z_k^{(1)}, \dots,
    \z_1^{(N)}, \dots, \z_k^{(N)} \big)
    \in \Zc^{k \times N}$}
    \For{$\sigma \in \Iknpr$} {
        Compute the estimator $\hat \theta \big(\z'_{\sigma(1)}, \dots, \z'_{\sigma(k)} \big)$ using the sample $(\X_i, \Z_i)$, $i=1,\dots,n$ \;
    }
    Compute the minimizer $\hat \beta$ of (\ref{condUstat:def:estimator_hat_beta}) using the $\hat \theta \big(\z'_{\sigma(1)}, \dots, \z'_{\sigma(k)} \big),$ $j=1,\dots,n',$ estimated in the above step \;
    \For{$i \leftarrow 1$ \KwTo $N$} {
         Compute the prediction
         $\tilde \theta(\z_1^{(i)}, \dots, \z_k^{(i)} )
         := \Lambda^{(-1)} \big( \psibm (\z_1^{(i)}, \dots, \z_k^{(i)} )^T \hat \beta \big)$ \;
    }
    \KwOut{An estimator $\hat \beta$ and $N$ predictions $\tilde \theta(\z_1^{(i)}, \dots, \z_k^{(i)} ),$ $i=1, \dots, N$.}
\caption{Two-step estimation of $\beta$ and prediction of the conditional parameters $\theta(\z_1^{(i)}, \dots, \z_k^{(i)} ),$ for $i=1, \dots, N$}
\end{algorithm}

\mds

Once an estimator $\hat \beta$ of $\beta^*$ has been computed, the prediction of all the conditional functionals is reduced to the computation of $\Lambda^{(-1)} \big( \psibm (\z_1^{(i)}, \dots, \z_k^{(i)} )^T \hat \beta \big) := \tilde \theta(\z_1^{(i)}, \dots, \z_k^{(i)} )$, for every $i=1, \dots, N$.
The total computational cost of this new method is therefore
$O(|\Iknpr|  n'{}^k + |\Iknpr| r + N s)$ operations. The first term corresponds to the cost of evaluating each non-parametric estimator (\ref{condUstat:def:estimator_hat_theta}). The second term corresponds to the minimization of the convex optimization program (\ref{condUstat:def:estimator_hat_beta}), and the last one is the prediction cost.
Note that the procedure described in Algorithm~\ref{condUstat:algo:estimation_beta} can provide a huge improvement compared to the previously available estimator with a cost in $O(N n^k)$ when $N \to \infty$, i.e. when we want to recover the full function $\theta(\cdot, \cdots, \cdot)$. Moreover, the speed-up given by Algorithm~\ref{condUstat:algo:estimation_beta} compared to the original conditional U-statistics (\ref{condUstat:def:estimator_hat_theta}) even increases with the sample size $n$, for moderate choices of $n'$.

A similar model, called \textit{functional response}, has already been studied: see, e.g. Kowalski and Tu~\cite[Chapter 6.2]{kowalski2008modern}. They provide a method to estimate the parameter $\beta^*$, using generalized estimating equations. However, they only provide asymptotic results for their estimator, and their algorithm needs to solve a multi-dimensional equation which has no reason to be convex.

\mds

In Section~\ref{section:theoretical_prop_theta}, we provide non-asymptotic bounds for the non-parametric estimator $\hat \theta$. Then Section~\ref{section:theoretical_prop_beta} is devoted to the statement of corresponding bounds, as well as asymptotic properties for the parametric estimator $\hat \beta$. Finally, a few examples are presented in Section~\ref{condUstats:sec:applications}. All proofs have been postponed to the Appendix.

\section{Theoretical properties of the nonparametric estimator \texorpdfstring{$\hat \theta(\cdot)$}{hat theta(.)}}
\label{section:theoretical_prop_theta}


%
\subsection{Non-asymptotic bounds for \texorpdfstring{$N_k$}{Nk}}

We remark that the estimator $\hat \theta$ is well-defined if and only if $N_k(\z_1, \dots, \z_k) > 0$, where
\begin{equation}
    N_k(\z_1, \dots, \z_k) := \frac{k!(n-k)!}{n!} \sum_{\sigma \in \Iknup} K_h \big(\Z_{\sigma(1)}-\z_1 \big) \cdots K_h \big(\Z_{\sigma(k)}-\z_k \big).
    \label{condUstat:def:normalization_factor_N_k}
\end{equation}
To prove that our estimator $\hat \theta(\z_1, \dots, \z_k)$ exists with a probability that tends to 1, we will therefore study the behavior of $N_k$. We will need the following assumptions to control the kernel $K$ and the density of $\Z$.

\begin{assumpt}
    The kernel $K(\cdot)$ is bounded, i.e. there exists a finite constant $C_K$ such that $K(\cdot) \leq C_{K}$ and $\int K(\u) d\u = 1$. The kernel is of order $\alpha$ for some $\alpha > 0$,
    i.e. for all $j = 1, \dots, \alpha - 1$ and all $1 \leq i_1, \dots, i_{\alpha} \leq p$,
    $\int K(\u) \, u_{i_1} \dots u_{i_j} \; d\u = 0.$
    \label{condUstat:assumpt:kernel_integral}
\end{assumpt}

\begin{assumpt}
    $f_\Z$ is $\alpha$-times continuously differentiable on $\Zc$ and there exists a finite constant $C_{K,\alpha}$ such that, for all $\z_1, \dots \z_k$,
    \begin{align*}
        &\int
        \Big| K \big(\u_1 \big) \cdots K \big(\u_k \big) \Big| 
        \sum_{m_1 + \, \cdots \, + m_k = \alpha}
        \binom{\alpha}{m_{1:k}} \\
        &\hspace{1cm} \cdot \prod_{i=1}^k \sum_{j_1, \dots, j_{m_i} = 1}^{p}
        \big| u_{i,j_1} \dots u_{i,j_{m_i}} \big| 
        \sup_{t \in [0,1]} \left| \frac{\partial^{m_i} f_{\Z}}{\partial z_{j_1} \, \cdots \, \partial z_{j_{m_i}}}
        \big(\z_i + t \u_i \big) \right| d\u_1 \cdots d\u_k
        \leq C_{K,\alpha}
    \end{align*}
    where $\binom{\alpha}{m_{1:k}} := \alpha! / \big( \prod_{i=1}^k (m_i!) \big)$ is the multinomial coefficient.
    \label{condUstat:assumpt:bound_int_deriv_fZ}
\end{assumpt}

\begin{assumpt}
    $f_\Z(\cdot) \leq f_{\Z, max} $ for some finite constant $f_{\Z, max}$.
    \label{condUstat:assumpt:f_Z_max}
\end{assumpt}

\begin{lemma}
    Under Assumptions \ref{condUstat:assumpt:kernel_integral}, \ref{condUstat:assumpt:bound_int_deriv_fZ} and \ref{condUstat:assumpt:f_Z_max},
    we have for any $t > 0$,
    \begin{align*}
        \PP \bigg( \big| N_k(\z_1, \dots, \z_k) &- \prod_{i=1}^k f_{\Z}(\z_i) \big| \leq \frac{C_{K, \alpha}}{\alpha!} h^{\alpha} + t \bigg)
        \geq 1 - 2 \exp \bigg( - \frac{[ n / k] t^2}{h^{-k p} C_1 + h^{-k p} C_2 t} \bigg),
    \end{align*}
    where $C_1 := 2 f_{\Z, max}^k ||K||_2^{2k}$,
    and $C_2 := (4/3) C_K^k$ and $||K||_2^2 := \int K^2$.
    \label{lemma:bound_Nk_fk}
\end{lemma}

This Lemma is proved in~\ref{condUstat:proof:lemma:bound_Nk_fk}.
More can be said if the density $f_\Z$ is bounded below. Therefore, we will use the following assumption.

\begin{assumpt}
    There exists a constant $f_{\Z, min} > 0$ such that
    for every $\z \in \Zc$, $f_{\Z}(\z) > f_{\Z, min}$.
    \label{condUstat:assumpt:f_Z_min}
\end{assumpt}

If for some $\epsilon > 0$, we have
$C_{K, \alpha} h^{\alpha} / \alpha ! + t
\leq f_{\Z, min} - \epsilon$,
then $\hat f(\z) \geq \epsilon > 0$ with probability larger than on the event whose probability is bound in Lemma \ref{lemma:bound_Nk_fk}.
We should therefore choose the largest $t$ possible, which yields the following corollary.

\begin{cor}
    Under Assumptions \ref{condUstat:assumpt:kernel_integral}-\ref{condUstat:assumpt:f_Z_min}, if
    $C_{K, \alpha} h^{\alpha} / \alpha ! \, < f_{\Z, min}$,
    then the random variable $N_k(\z_1, \dots, \z_k)$ is strictly positive with a probability larger than
    $1 - 2 \exp \bigg( - \frac{[ n / k] h^{k p}  \big( f_{\Z, min} - C_{K, \alpha} h^{\alpha}  / \alpha ! \big)^2}
    {C_1 + C_2 \big( f_{\Z, min} - C_{K, \alpha} h^{\alpha}  / \alpha ! \big) } \bigg),$
    guaranteeing the existence of the estimator $\hat \theta(\z_1, \dots, \z_k)$ on this event.
    \label{cor:probaTheta_Z_valid}
\end{cor}

\subsection{Non-asymptotic bounds in probability for \texorpdfstring{$\hat \theta$}{hat theta}}

In this section, we generalize the bounds given in~\cite{derumigny2018kernel} for the conditional Kendall's tau to any conditional U-statistics.
To establish bounds on $\hat \theta$ for every fixed $n$, we will need some assumptions on the joint law of $(\X, \Z)$.
\begin{assumpt}
    There exists a measure $\mu$ on $(\Xc, \Ac)$ such that $\PP_{\X,\Z}$ is absolutely continuous with respect to $\mu \otimes Leb_p$, where $Leb_p$ is the Lebesgue measure on $\Rb^p$.
    \label{condUstat:assumpt:dominating_measure_mu}
\end{assumpt}

\begin{assumpt}
    For every $\x \in \Xc$, $\z \mapsto f_{\X, \Z}(\x, \z)$ is differentiable almost everywhere up to the order $\alpha$.
    Moreover, there exists a finite constant $C_{g,f, \alpha} > 0$, such that, for every positive integers $m_1, \dots, m_k$ such that $\sum_{i=1}^k m_i = \alpha$, for every $0 \leq j_1, \dots, j_{m_i} \leq p$,
    \begin{align*}
        &\int \prod_{i=1}^k \Bigg|
        \Big( g \big(\x_1, \dots, \x_k \big)
        - \EE \big[ g(\X_1, \dots, \X_k) \big| \Z_i=\z_i, \forall i=1, \dots, k \big] \Big) \\
        &\hspace{0cm} \cdot \left( \frac{\partial^{m_i} f_{\X,\Z}}{\partial z_{j_1} \, \cdots \, \partial z_{j_{m_i}}} \big(\x_i,\z_i + \u_i \big) - \frac{\partial^{m_i} f_{\X,\Z}}{\partial z_{j_1} \, \cdots \, \partial z_{j_{m_i}}} \big(\x_i,\z_i \big) \right) \Bigg| d\mu(\x_1) \cdots d\mu(\x_k)
        \leq C_{g, f, \alpha} \prod_{i=1}^k \big|\u_i \big|_\infty,
    \end{align*}
    for every choices of $\x_1, \dots, \x_k \in \Xc$ and $\z_1, \dots, \z_k \in \Zc, \u_1, \dots, \u_k \in \Rb^p$ such that $\z_i + \u_i \in \Zc$. There exists a constant $C'_{K,\alpha}$ such that 
    ${\sum_{m_1 + \, \cdots \, + m_k = \alpha}
    \binom{n}{m_{1:k}}
    \int \prod_{i=1}^k K(\u_i)
    \sum_{j_1, \dots, j_{m_i} = 1}^{p} u_{i,j_1} \dots u_{i,j_{m_i}} \prod_{i=1}^k \big|\u_i \big|_\infty d\u_1 \cdots
    d\u_k \leq C'_{K, \alpha}.}$
    \label{condUstat:assumpt:f_XZ_Holder}
\end{assumpt}
\vspace{-0.5cm}
An easy situation is the case when $g$ is bounded, i.e. when the following assumption hold.
\begin{assumpt}
    There exists a constant $C_{g}$ such that $||g||_\infty \leq C_{g} < + \infty$.
    \label{condUstat:assumpt:g_bounded}
\end{assumpt}
When $g$ is not bounded, a weaker result can still be proved under a ``conditional Bernstein'' assumption. This assumption will help us to control the tail behavior of $g$ so that exponential concentration bounds are available.
\begin{assumpt}[conditional Bernstein assumption]
    There exists a positive function $B_{g}$ such that, for all $l \geq 1$ and $\z_1, \dots, \z_k \in \Rb^{k p}$,
    $\EE \Big[ \big| g(\X_1, \dots, \X_k) \big|^l
        \, \Big| \Z_1 = \z_1, \dots, \Z_k = \z_k \Big]
        \leq B_{g}(\z_1, \dots, \z_k)^l l!,$
    such that $B_{g}(\Z_1, \dots, \Z_k) \leq \tilde B_{g}$ almost surely, for some finite positive constant $\tilde B_{g}$.
    \label{condUstat:assumpt:conditional_Bernstein}
\end{assumpt}
As a shortcut notation, we will define also $B_{g, \z} := B_{g}(\z_1, \dots, \z_k)$. The following proposition is proved in~\ref{condUstat:proof:prop:exponential_bound_hat_theta}.

\begin{prop}[Exponential bound for the estimator $\hat \theta(\z_1, \dots, \z_k)$, with fixed $\z_1, \dots \z_k \in \Zc^k$]
    Assume either Assumption \ref{condUstat:assumpt:g_bounded} or the weaker Assumption \ref{condUstat:assumpt:conditional_Bernstein}.
    Under Assumptions \ref{condUstat:assumpt:kernel_integral}-\ref{condUstat:assumpt:f_XZ_Holder},
    for every $t,t'>0$ such that
    $C_{K, \alpha} h^{\alpha} / \alpha ! + t < f_{\Z, min}/2$, we have
    \begin{align*}
        &\PP \Bigg( \big| \hat \theta(\z_1, \dots, \z_k) - \theta(\z_1, \dots, \z_k) \big|
        < \big(1 + C_3 h^{\alpha} + C_4 t \big) \times
        \big( C_5 h^{k + \alpha} + t' \big) \Bigg) \\
        & \hspace{2cm} \geq 1
        - 2 \exp \bigg( - \frac{[ n / k] t^2 h^{k p}}{ C_1 + C_2 t} \bigg)
        - 2 \exp \bigg( - \frac{[ n / k] t'{}^2 h^{kp}}
        {C_6 + C_7 t'} \bigg),
    \end{align*}
    where
    $C_3 := 4 f_{\Z, max}^k f_{\Z, min}^{-2k} C_{K, \alpha} / \alpha !$,
    $C_4 := 4 f_{\Z, max}^k f_{\Z, min}^{-2k} $ and
    $C_5 := C_{g, f, \alpha} C'_{K,\alpha}
    f_{\Z,min}^{-k} / \alpha !$.
    
    \mds
    
    If Assumption \ref{condUstat:assumpt:g_bounded} is satisfied, the result holds with the following values:
    $C_6 := 2 C_{g}^2 f_{\Z,max}^k f_{\Z,min}^{-2k} ||K||_2^{2k}$,
    $C_7 := (8/3) C_K^k C_{g}^k f_{\Z,min}^{-k}$ ;
    in the case of Assumption \ref{condUstat:assumpt:conditional_Bernstein}, the result holds with the following alternative values:
    $\tilde C_6 := 128 \big(B_{g, \z} + \tilde B_{g} \big)^2 C_K^{2k-1} f_{\Z,min}^{-2k}$,
    $\tilde C_7 := 2 \big(B_{g, \z} + \tilde B_{g} \big) C_K^{k} f_{\Z,min}^{-k}$.
    \label{condUstat:prop:exponential_bound_hat_theta}    
\end{prop}

\section{Theoretical properties of the estimator \texorpdfstring{$\hat \beta$}{hat beta}}
\label{section:theoretical_prop_beta}



Let us define the matrix $\Zb$ of dimension $|\Iknpr| \times r$ by $[\Zb']_{i,j} := \psi_j \big(\z'_{\sigma_i(1)}, \dots, \z'_{\sigma_i(k)} \big)$, where $1 \leq i \leq |\Iknpr|$, $1 \leq j \leq r$ and $\sigma_i$ is the $i$-th element of $\Iknpr$. The chosen order of $\Iknpr$ is arbitrary and has no impact in practice. In the same way, we define the vector $\Y$ of dimension $|\Iknpr|$ defined by
$Y_i := \Lambda \Big(\hat \theta \big(\z'_{\sigma_i(1)}, \dots, \z'_{\sigma_i(k)} \big) \Big),$ such that the criterion (\ref{condUstat:def:estimator_hat_beta}) is in the standard Lasso form $\hat \beta := \arg \min_{\beta \in \Rb^{r}}
\Big[ ||\Y - \Zb' \beta||^2 + \lambda |\beta|_1 \Big],$
where for any vector $\v$ of size $|\Iknpr|$, its scaled norm is defined by $||\v|| := |\v|_2 / \sqrt{|\Iknpr|}$.
Following~\cite{derumigny2018kendall}, we define $\xi_{i,n}$, for $1 \leq i \leq |\Iknpr|$, by
$\xi_{i,n} = \xi_{\sigma_i,n}
:= \Lambda \Big( \hat \theta
\big(\z'_{\sigma_i(1)}, \dots, \z'_{\sigma_i(k)} \big) \Big)
- \psibm \big(\z'_{\sigma_i(1)}, \dots,
\z'_{\sigma_i(k)} \big)^T \beta^*.$

\subsection{Non-asymptotic bounds on \texorpdfstring{$\hat \beta$}{hat beta}}

We will also use the \emph{Restricted Eigenvalue} (RE) condition, introduced by Bickel, Ritov and Tsybakov~\cite{bickel2009simultaneous}.
For $c_0 > 0$ and $s\in \{ 1, \dots, p \}$, it is defined as follows:

\medskip

\noindent
$RE(s,c_0)$ \textbf{condition :}
\emph{The design matrix $\Zb'$ satisfies}
$$\kappa(s, c_0) := \min
\left\{
\frac{||\Zb' \delta||}{|\delta|_2}
: \delta \neq 0, \, |\delta_{J_0^C}|_1 \leq c_0 |\delta_{J_0}|_1 , \, 
J_0 \subset \{1, \dots, r\}, \, |J_0| \leq s
\right\} > 0.$$

Note that this condition is very mild, and is satisfied with a high probability for a large class of random matrices: see Bellec et al. \cite[Section 8.1]{bellec2018slope} for references and a discussion.
We will also need the following regularity assumption on the function $\Lambda(\cdot)$.
\begin{assumpt}
    The function $\z\mapsto \psibm(\z)$ are bounded on $\Zc$ by a constant $C_\psibm$. Moreover, $\Lambda(\cdot)$ is continuously differentiable.
    Let $\Tc$ be the range of $\theta$, from $\Zc^k$ towards $\Rb$.
    On an open neighborhood of $\Tc$, the derivative of $\Lambda(\cdot)$ is bounded by a constant $C_{\Lambda'}$.
    \label{condUstat:assumpt:compact_case}
\end{assumpt}

The following theorem is proved in~\ref{condUstat:proof:thm:bound_proba_hat_beta}.

\begin{thm}
    Assume either Assumption \ref{condUstat:assumpt:g_bounded} or the weaker Assumption \ref{condUstat:assumpt:conditional_Bernstein}.
    Suppose that Assumptions \ref{condUstat:assumpt:kernel_integral}-\ref{condUstat:assumpt:f_XZ_Holder} and \ref{condUstat:assumpt:compact_case}
    hold and that the design matrix $\Zb'$ satisfies the $RE(s,3)$ condition.
    Choose the tuning parameter as $\lambda = \gamma t$,
    with $\gamma \geq 4$ and $t>0$, and assume that we choose $h$ small enough such that
    \begin{align}
        h &\leq
        \min \bigg( \Big( \frac{f_{\Z, min} \alpha !}
        {4 \,  C_{K, \alpha}} \Big)^{1/\alpha} \; , \;
        \Big( \frac{t}{2 C_5 C_8} \Big)^{1/(k+\alpha)}
        \bigg)
        \label{condUstat:cond:h_1_2},
    \end{align}
    where $C_8 := C_{\psibm} C_{\Lambda'}
    \big(1 + C_4 f_{\Z, min} / 2 \big)$.
    Then, we have
    \begin{align}
        &\PP \Big( ||\Zb'(\hat \beta - \beta^*)||
        \leq \dfrac{4(\gamma+1) t \sqrt{s}}{\kappa(s,3)}
        \text{ and } |\hat \beta - \beta^*|_q
        \leq \dfrac{4^{2/q}(\gamma+1) t s^{1/q} }{\kappa^2(s,3)},
        \text{ for every } 1 \leq q \leq 2 \Big) \nonumber \\
        &\hspace{1cm} \geq 1 -  2 \sum_{\sigma \in \Iknpr} \Bigg[
        \exp \bigg( - \frac{[ n / k] f_{\Z, min}^2 h^{k p}}
        { 16 C_1 + 4 C_2 f_{\Z, min} } \bigg)
        + \exp \bigg( - \frac{[ n / k] t^2 h^{kp}}
        {4 C_8^2 C_{6,\sigma} + 2 C_8 C_{7,\sigma} t} \bigg) \Bigg].
        \label{condUstat:eq:bound_proba_hat_beta}
    \end{align}
    If Assumption \ref{condUstat:assumpt:g_bounded} is satisfied, the result holds with
    $C_{6,\sigma}$ and $C_{7,\sigma}$ constant, respectively to $C_6$ and $C_7$ defined in Proposition \ref{condUstat:prop:exponential_bound_hat_theta}.
    In the case of Assumption \ref{condUstat:assumpt:conditional_Bernstein}, the result holds with the following alternative values:
    $C_{6,\sigma} := 128 \big(B_{g}(\z'_{\sigma(1)}, \dots, \z'_{\sigma(k)}) + \tilde B_{g} \big)^2 C_K^{2k} f_{\Z,min}^{-2k}$ and
    $C_{7,\sigma} := 2 \big(B_{g}(\z'_{\sigma(1)}, \dots, \z'_{\sigma(k)}) + \tilde B_{g} \big) C_K^{k} f_{\Z,min}^{-k}$.
    \label{condUstat:thm:bound_proba_hat_beta}
\end{thm}

The latter theorem gives some bounds that hold in probability for the prediction error $||\Zb'(\hat \beta - \beta^*)||_{n'}$ and for the estimation error $|\hat \beta - \beta^*|_q$ with $1 \leq q \leq 2$ under the specification (\ref{model:lambda_cond_theta_Z}). Note that the influence of $n'$ and $r$ is hidden through the Restricted Eigenvalue number $\kappa(s,3)$.

\subsection{Asymptotic properties of \texorpdfstring{$\hat \beta$}{hat beta} when \texorpdfstring{$n\to\infty$}{n to infty} and for fixed \texorpdfstring{$n'$}{n prime}}
\label{sec:Asymp_n_fixed_nprime}

In this part, $n'$ is still assumed to be fixed and we state the consistency and the asymptotic normality of $\hat\beta$ as $n \to \infty$.
As above, we adopt a fixed design: the $\z'_i$ are arbitrarily fixed or, equivalently, our reasoning are made conditionally on the second sample. In this section, we follow Section 3 of Derumigny and Fermanian~\cite{derumigny2018kendall} which gives similar results for the conditional Kendall's tau, a particular conditional U-statistic of order~$2$. Proofs are identical and therefore omitted.
Nevertheless, asymptotic properties of $\hat \beta$ require corresponding results on the first-step estimators $\hat \theta$. These results are state in Stute~\cite{stute1991conditional} and recalled for convenience in~\ref{section:asymptotic_result_theta_stute}.
For $n, n' > 0$, denote by $\hat \beta_{n, n'}$ the estimator (\ref{condUstat:def:estimator_hat_beta}) with $h = h_{n}$ and $\lambda = \lambda_{n, n'}$.

\begin{lemma}
    We have $\hat \beta_{n, n'} = \arg \min_{\beta \in \Rb^{p'}}
    \GG_{n,n'}(\beta)$, where
    \begin{align}
        \GG_{n,n'}(\beta)
        &:= \frac{2 (n'-k)!}{n'!} \sum_{\sigma \in \Iknpr} \xi_{\sigma,n} \psibm \big(\z'_{\sigma(1)}, \dots, \z'_{\sigma(k)} \big)^T (\beta^* - \beta) \nonumber \\
        &+ \frac{(n'-k)!}{n'!} \sum_{\sigma \in \Iknpr}
        \big\{ \psibm \big(\z'_{\sigma(1)}, \dots, \z'_{\sigma(k)} \big)^T (\beta^* - \beta) \big\}^2
        + \lambda_{n, n'} |\beta|_1.
        \label{condUstat:def:process_GG}
    \end{align}
    \label{lemma:beta_process_GG}
\end{lemma}
\vspace{-1cm}
\begin{thm}[Consistency of $\hat \beta$]
    Under Assumption~\ref{condUstat:assumpt:stute_consistency},
    if $n'$ is fixed and $\lambda = \lambda_{n, n'} \to \lambda_0 $,
    then, given $\z'_1, \ldots, \z'_{n'}$ and as $n$ tends to the infinity, $\hat \beta_{n, n'} \inprobto \beta^{**} := \inf_\beta \GG_{\infty,n'} (\beta),$ where
    \begin{align*}
        \GG_{\infty,n'} (\beta)
        := \frac{1}{n'} \sum_{\sigma \in \Iknpr}
        \Big( \psibm \big(\z'_{\sigma(1)}, \dots, \z'_{\sigma(k)} \big)^T (\beta^* - \beta) \Big)^2 + \lambda_0 |\beta |_1.
    \end{align*}
    In particular, if $\lambda_0=0$ and
    $< \{ \psibm \big(\z'_{\sigma(1)}, \dots, \z'_{\sigma(k)} \big) : \sigma \in \Iknpr \} > \, = \Rb^{r}$, then
    $\hat \beta_{n, n'} \inprobto \beta^*$.
\end{thm}

\begin{thm}[Asymptotic law of the estimator]
\label{condUstat:thm:WeakConvLasso}
    Under Assumption~\ref{condUstat:assumpt:stute_asymptotic_normality}, and if
    $\lambda_{n, n'}  (n h_{n,n'}^{p})^{1/2}$ tends to $ \ell $ when $n\rightarrow \infty$, we have
    $(n h_{n,n'}^{p})^{1/2} (\hat \beta_{n, n'} - \beta^*)
    \indistrto \u^* := \arg \min_{\u \in \Rb^{r}} \FF_{\infty, n'}(\u),$
    given $\z'_1,\ldots, \z'_{n'}$, where
    \begin{align*}
        \FF_{\infty, n'}(\u)
        &:= \frac{2 (n'-k)!}{n'!} \sum_{\sigma \in \Iknpr} \sum_{j=1}^{r}
        W_{\sigma} \psi_j \big(\z'_{\sigma(1)}, \dots, \z'_{\sigma(k)} \big) u_{j}
        + \frac{(n'-k)!}{n'!} \sum_{\sigma \in \Iknpr} \left( \psibm \big(\z'_{\sigma(1)}, \dots, \z'_{\sigma(k)} \big)^T \u \right)^2 \\
        &+ \ell \sum_{i=1}^{r} \big( |u_i| \1_{\{\beta_i^*=0\}} + u_i\sgn(\beta_i^*) \1_{\{\beta_i^* \neq 0\}} \big) ,
    \end{align*}
    with $\W = (W_\sigma)_{\sigma \in \Iknpr} \sim \Nc \left(0, \tilde \HH \right)$ where
    \begin{align*}
        [\tilde \HH]_{\sigma, \varsigma}
        & := \sum_{j,l = 1}^k
        \1_{\left\{\z'_{\sigma(j)} = \z'_{\varsigma(l)}\right\}} 
        \frac{||K||_2^2}{f_\Z \left( \z'_{\sigma(j)} \right)} 
        \Lambda' \bigg(\theta \left(\z'_{\sigma(1)}, \dots, \z'_{\sigma(k)} \right) \bigg)
        \Lambda' \bigg(\theta \left(\z'_{\varsigma(1)}, \dots, \z'_{\varsigma(k)} \right) \bigg) \\
        & \cdot \bigg(\tilde \theta_{j,l}
        \left(\z'_{\sigma(1)}, \dots, \z'_{\sigma(k)} ,
        \z'_{\varsigma(1)}, \dots, \z'_{\varsigma(k)} \right)
        - \theta \left(\z'_{\sigma(1)}, \dots, \z'_{\sigma(k)} \right)
        \theta \left(\z'_{\varsigma(1)}, \dots, \z'_{\varsigma(k)} \right) \bigg),
    \end{align*}
    and $\tilde \theta_{j,l}$ is as defined in Equation (\ref{condUstat:eq:def:tilde_theta_jl}).
    
    \mds
    
    Moreover, $\lim\sup_{\protect n \to \infty} \PP\left( \Sc_n= \Sc \right) = c < 1, $
    where $\Sc_n := \{j:\hat\beta_j \neq 0\}$ and $\Sc := \{j:\beta_j \neq 0\}$.
\end{thm}

\medskip

A usual way of obtaining the oracle property is to modify our estimator in an ``adaptive'' way. Following Zou~\cite{zou2006}, consider a preliminary ``rough'' estimator of $\beta^*$, denoted by $\tilde \beta_n$, or more simply $\tilde\beta$. Moreover $\nu_n (\tilde \beta_n - \beta^*)$ is assumed to be asymptotically normal, for some deterministic sequence $(\nu_n)$ that tends to the infinity.
Now, let us consider the same optimization program as in~(\ref{condUstat:def:estimator_hat_beta}) but with a random tuning parameter given by
$   \lambda_{n,n'} := \tilde\lambda_{n,n'}/ |\tilde{\beta}_n|^\delta$,
for some constant $\delta >0$ and some positive deterministic sequence $(\tilde\lambda_{n,n'})$. The corresponding adaptive estimator (solution of the modified Equation~(\ref{condUstat:def:estimator_hat_beta})) will be denoted by $\check{\beta}_{n,n'}$, or simply $\check{\beta}$. Hereafter, we still set
$\Sc_n=\{j : \check{\beta}_j \neq 0\}$.

\begin{thm}[Asymptotic law of the adaptive estimator of $\beta$]
    Under Assumption~\ref{condUstat:assumpt:stute_asymptotic_normality}, if 
    $\tilde\lambda_{n, n'} (n h_{n,n'}^{p})^{1/2}  \rightarrow \ell\geq 0 $ and
    $\tilde\lambda_{n, n'} (n h_{n,n'}^p)^{1/2} \nu_n^\delta  \rightarrow \infty $ when $n\rightarrow \infty$, we have $(n h_{n,n'}^{p})^{1/2} (\check \beta_{n, n'} - \beta^*)_{\Sc}
    \indistrto \u_{\Sc}^{**} :=
    \underset{\u_{\Sc} \in \Rb^{s}}{\arg \min} \,
    \check\FF_{\infty, n'}(\u_{\Sc}),$ where
    \begin{align*}
        &\check\FF_{\infty, n'}(\u_{\Sc})
        := \frac{2 (n'-k)!}{n'!} \sum_{\sigma \in \Iknpr}
        \sum_{j\in \Sc} W_\sigma \psi_j (\z'_i) u_{j}
        + \frac{(n'-k)!}{n'!} \sum_{\sigma \in \Iknpr} \Big(
        \sum_{j\in \Sc} \psi_j (\z'_i) u_{j} \Big)^2
        + \ell \sum_{i\in \Sc} \frac{u_i}{|\beta_i^*|^\delta} \sgn(\beta_i^*) ,
    \end{align*}
    and $\W = (W_\sigma)_{\sigma \in \Iknpr} \sim \Nc \big(0, \tilde \HH \big).$
    
    \noindent
    Moreover, when $\ell=0$, the oracle property is fulfilled: $ \PP\left( \Sc_n= \Sc \right) \rightarrow 1$ as $n \to \infty$.
    \label{condUstat:thm:WeakConvLassoAdaptive}
\end{thm}

\subsection{Asymptotic properties of \texorpdfstring{$\hat \beta$}{hat beta} jointly in \texorpdfstring{$(n,n')$}{n,n prime}}

Now, we consider the framework in which both $n$ and $n'$ are going to infinity, while the dimensions $p$ and~$r$ stay fixed.
We now provide a consistency result for $\hat \beta_{n,n'}$.
\begin{thm}[Consistency of $\hat \beta_{n,n'}$, jointly in $(n,n')$]
    Assume that Assumptions \ref{condUstat:assumpt:kernel_integral}-\ref{condUstat:assumpt:f_XZ_Holder}, \ref{condUstat:assumpt:conditional_Bernstein} and \ref{condUstat:assumpt:compact_case} are satisfied.
    Assume that
    $\sum_{\sigma \in \Iknpr}
    \psibm \big(\z'_{\sigma(1)}, \dots, \z'_{\sigma(k)} \big)
    \psibm \big(\z'_{\sigma(1)}, \dots, \z'_{\sigma(k)} \big)^T/n'$ converges to a matrix $M_{\psi,\z'}$, as $n' \to \infty$.
    Assume that $\lambda_{n, n'} \to \lambda_{0}$
    and $n' \exp ( - A n h^{2 k p}) \to 0$ for every $A>0$, when $(n,n') \to \infty$.
    Then $\hat \beta_{n,n'} \inprobto \arg \min_{\beta \in \Rb^{r}} \GG_{\infty,\infty}(\beta),$ as $(n,n') \to \infty,$
    where $\GG_{\infty,\infty}(\beta):= (\beta^* - \beta) M_{\psi,\z'} (\beta^* - \beta)^T + \lambda_{0} |\beta|_1$.
    Moreover, if $\lambda_0 = 0$ and $M_{\psi,\z'}$ is invertible, then $\hat \beta_{n,n'}$ is consistent and tends to the true value~$\beta^*$.
    \label{condUstat:thm:consistency_hatBeta_n_nprime}
\end{thm}

Note that, since the sequence $(\z'_i)$ is deterministic, we only assume the convergence of the sequence of deterministic matrices
$\sum_{\sigma \in \Iknpr}
\psibm \big(\z'_{\sigma(1)}, \dots, \z'_{\sigma(k)} \big)
\psibm \big(\z'_{\sigma(1)}, \dots, \z'_{\sigma(k)} \big)^T/n'$
in $\Rb^{r^2}$.
Moreover, if the ``second subset'' $(\z'_i)_{i=1,\ldots,n'}$ were a random sample (drawn along the law $\PP_{\Z}$), the latter convergence would be understood ``in probability''.
And if $\PP_{\Z}$ satisfies the identifiability condition (Proposition \ref{condUstat:prop:identifiability_condition}), then $M_{\psi,\z'}$ would be invertible and $\hat \beta_{n,n'}\to \beta^*$ in probability.
Now, we want to go one step further and derive the asymptotic law of the estimator $\hat \beta_{n,n'}$.

\begin{thm}[Asymptotic law of $\hat \beta_{n,n'}$, jointly in $(n,n')$]
    Under Assumptions \ref{condUstat:assumpt:kernel_integral}-\ref{condUstat:assumpt:dominating_measure_mu}
    and under Assumption~\ref{condUstat:assumpt:asymptNorm_joint},
    we have
    $$ (n \times n' \times h_{n,n'}^{p})^{1/2} (\hat \beta_{n, n'} - \beta^*) \indistrto \Nc(0, \tilde V_{as}),$$
    where $\tilde V_{as} := V_1^{-1} V_2 V_1^{-1}$, $V_1$ is the matrix defined in Assumption~\ref{condUstat:assumpt:asymptNorm_joint}(iv), and $V_2$ in Assumption~\ref{condUstat:assumpt:asymptNorm_joint}(v).
    \label{condUstat:thm:weak_conv_doubleAsympt}
\end{thm}

This theorem is proved in~\ref{condUstat:proof:thm:weak_conv_doubleAsympt} where we state Assumption~\ref{condUstat:assumpt:asymptNorm_joint}.

\section{Applications and examples}
\label{condUstats:sec:applications}


Following Example 4.4 in Stute~\cite{stute1991conditional}, we consider the function
$g(x_1, x_2) := \1 \{ x_1 \leq x_2 \}$, with $k=2$.
In this case $\theta(\z_1, \z_2) = \PP( X_1 \leq X_2 | \Z_1 = \z_1, \Z_2 = \z_2)$. The parameter $\theta(\z_1, \z_2)$ quantifies the probability that the quantity of interest $X$ be smaller if we knew that $\Z = \z_1$ than if we knew that $\Z = \z_2$.

\mds

To illustrate our methods, we choose a simple example, with the Epanechnikov kernel, defined by $K(x) := (3/4)(1-u^2)\1{|u| \leq 1}$. It is a kernel of order $\alpha = 2$, with $\int K^2 = 3/5$. Assumption~\ref{condUstat:assumpt:kernel_integral} is then satisfied with $C_K := 3/4$.
Fix $p=1$, $\Zc = [-1,1]$, $\Xc = \Rb$,
$f_Z(z) = \phi(z) \1 \{|z| \leq 1 \} / (1 - 2 \Phi(-1))$, where $\Phi$ and $\phi$ are respectively the cdf and the density of the standard Gaussian distribution and $X|Z=z \sim \Nc(z, 1)$, for every $z \in \Zc$.

Assumption~\ref{condUstat:assumpt:bound_int_deriv_fZ} is then satisfied with $C_{K,\alpha} = 0.2
$. Assumption~\ref{condUstat:assumpt:f_Z_max} is easily satisfied with $f_{Z,\max} = 1 / \big( \sqrt{2 \pi} (1-2\Phi(-1)) \big)
\leq 0.59$.
Therefore, we can apply Lemma~\ref{lemma:bound_Nk_fk}.
We compute the constants
$C_1 := 2 f_{\Z, max}^k ||K||_2^{2k}
= 2 \times 0.59^2 \times (3/5)^2
\leq 0.26$
and $C_2 := (4/3) C_K^k = (4/3) \times (3/4)^2 = 3/4$.
Therefore, for any $n \geq 0$, $h, t > 0$, $z_1, z_2 \in \Zc$, we have
\begin{align*}
    \PP \bigg( \big| N_2(z_1, z_2) &- f_{Z}(z_1) f_{Z}(z_2) \big|
    \leq 0.1 h^{\alpha} + t \bigg)
    \geq 1 - 2 \exp \bigg( - \frac{[ n / 2] t^2}
    {0.26 h^2 + 0.75 h^2 t} \bigg),
\end{align*}
Assumption~\ref{condUstat:assumpt:f_Z_min} is satisfied with
$f_{Z, \min} = \phi(1) / (1 - 2 \Phi(-1)) > 0.35$, so that we can apply Corollary~\ref{cor:probaTheta_Z_valid}.
Therefore, the estimator $\hat \theta(z_1, z_2)$ exists with probability greater than 
$1 - 2 \exp \bigg( - \frac{(n-1) h^2  \big( 0.35 - 0.1 h^2 \big)^2}
    {0.52 + 1.5 \times \big( 0.35 - 0.1 h^2 \big) } \bigg).$
Note that this probability is greater than $0.99$ as soon as
$n \geq 3 \big( 0.52 + 1.5 \times ( 0.35 - 0.1 h^2 ) \big)
/ \big(h^2 ( 0.35 - 0.1 h^2 )^2 \big)$. For example, with $h=0.2$, it means that the estimator $\hat \theta(z_1, z_2)$ exists with a probability greater than $99\%$ as soon as $n$ is greater than $651$.

\mds

We list below other possible examples of applications.
Conditional moments constitute also a natural class of U-statistics. They include the conditional variance ($p_{\X} = 1$, $k=2$, $g(X_1, X_2)=X_1^2 - X_1 \cdot X_2$) and the conditional covariance ($p_{\X} = 2$, $k=2$,
$g(\X_1, \X_2) := X_{1,1} \times X_{2,1} - X_{1,1} \times X_{2,2}$).
The conditional variance gives information about the volatility of $X$ given the variable $\Z$. Conditional covariances can be used to describe how the dependence moves as a function of the conditioning variables $\Z$.
Higher-order conditional moments (skewness, kurtosis, and so on) can also be estimated by higher-order conditional U-statistics, and they described respectively how the asymmetry and the behavior of the tails of $X$ change as function of $Z$. 

\mds

Gini's mean difference, an indicator of dispersion, can also be used in this framework. Formally, it is defined as the U-statistic with $p_{\X} = 1$, $k=2$ and $g(X_1, X_2) := |X_1 - X_2|$. Its conditional version describes how two variables are far away, on average, given their conditioning variables $\Z$. for example, $X$ could be the income of an individual, $\Z$ could be the position of its home, and $\theta(\z_1, \z_2)$ represent the average inequality between the income of two persons, one at point $\z_1$ and the other at point $\z_2$.

\mds

Other conditional dependence measures can also be written as conditional U-statistics, see e.g. Example 1.1.7 of Koroljuk and Borovskisch~\cite{koroljuk1994theory}. They show how a U-statistic of order $k=5$ can be used to estimated the dependence parameter $$\theta = \iint \big(F_{1,2}(x,y) - F_{1,2}(x,\infty) F_{1,2}(\infty,y)\big) dF_{1,2}(x,y).$$
In our framework, we could consider a conditional version, given by
$$\theta(\z_1, \z_2)
= \iint \big(F_{1,2|\Z=\z}(x,y) - F_{1,2|\Z=\z}(x,\infty) F_{1,2|\Z=\z}(\infty,y)\big) dF_{1,2|\Z=\z}(x,y),$$
where $\X$ is of dimension $p_{\X} = 2$.


\mds

\textbf{Acknowledgements: }
This work is supported by the GENES and by the Labex Ecodec under the grant ANR-11-LABEX-0047 from the French Agence Nationale de la Recherche. The author thanks Professor Jean-David Fermanian for helpful comments and discussions.


\bigskip

\bibliographystyle{abbrv}
\bibliography{biblio_Cond_UstatRegression}{}

\bigskip

\appendix

\section{Notations}

In the proofs, we will use the following shortcut notation.
First, $\x_{1:k}$ denotes the $k$-tuple ${(\x_1, \dots, \x_k) \in \Xc^k}$. Similarly, for a function $\sigma$, $\sigma(1:k)$ denotes the tuple $(\sigma(1), \dots, \sigma(k))$, and $\X_{\sigma(1:k)}$ is the $k$-tuple $(\X_{\sigma(1)}, \dots, \X_{\sigma(k)})$.
For any variable $Y$ and any collection of given points $(\z_1, \dots \z_k)$, the conditional expectation $\EE[Y | \Z_{1:k} = \z_{1:k}]$ denotes $\EE[Y | \Z_1 = \z_1, \dots, \Z_k = \z_k]$.
We denote by $\int \phi(\z_{1:k}) d\z_{1:k}$ the integral $\int \phi(\z_1, \dots, \z_k) d\z_1 \cdots d\z_k$ for any integrable function $\phi: \Rb^{k \times p} \to \Rb$, and by $\int g(\x_{1:k}) d\mu^{\otimes k}(\x_{1:k})$ the integral $\int g(\z_1, \dots, \z_k) d\mu(\x_1) \cdots d\mu(\x_k)$ for any $\mu$-integrable function $g: \Xc^k \to \Rb$.

\section{Asymptotic results for \texorpdfstring{$\hat \theta$}{hat theta}}
\label{section:asymptotic_result_theta_stute}

The estimator $\hat \theta(\z_1, \dots, \z_k)$ has been first studied by Stute (1991) \cite{stute1991conditional}. He proved the consistency and the asymptotic normality of $\hat \theta(\z_1, \dots, \z_k)$. We recall his results.

\begin{assumpt}
    \begin{enumerate}[(i)]
        \item $h_n \to 0$ and $n h_n^{p} \to \infty$ ;
        \item $K(\z) \geq C_{K,1} \1_{\{ |\z|_\infty \leq \, C_{K,2}\}}$ for some $C_{K,1}, \, C_{K,2} > 0$ ;
        \item there exists a decreasing function $H:\Rb_+ \to \Rb_+$, and positive constants $c_1, c_2$ such that
        $H(t) {\underset{t \to \infty}{=}} o(t^{-1})$ and $c_1 H(|\z|_\infty) \leq K(\z) \leq c_2 H(|\z|_\infty)$. 
    \end{enumerate}
    \label{condUstat:assumpt:stute_consistency}
\end{assumpt}

\begin{prop}[Consistency of $\hat \theta$, Theorem 2 in Stute \cite{stute1991conditional}]
    Under Assumption~\ref{condUstat:assumpt:stute_consistency}, for $\PP_\Z^{\otimes k}$-almost all $(\z_1, \dots, \z_k)$, $\hat \theta(\z_1, \dots, \z_k) \inprobto \theta(\z_1, \dots, \z_k)$ as $n \to \infty$.
    \label{condUstat:prop:consistency_hat_Tau}
\end{prop}

We introduce now a few more notations to state the asymptotic normality of $\hat \theta$.
For $1 \leq j,l,m \leq k$ and $\z_1, \dots, \z_{3k} \in \Zc^{3k}$, define
\begin{align}
    \theta_{j,l}&(\z_1, \dots, \z_k) :=
    \EE \big[
    g(\X_1, \dots, \X_{j-1}, \X, \X_{j+1}, \dots, \X_k)
    g(\X_{k+1}, \dots, \X_{k+l-1}, \X, \X_{k+l+1}, \dots, \X_{2k}) \nonumber \\
    &\big| \Z = \z_j \; ; \; \Z_i=\z_i, \forall i=1, \dots, k, i \neq j
    \; ; \; \Z_{k+i}=\z_i, \forall i=1, \dots, k, i \neq l \big], 
    \nonumber \\
    \tilde \theta_{j,l}&(\z_1, \dots, \z_{2k}) :=
    \EE \big[
    g(\X_1, \dots, \X_{j-1}, \X, \X_{j+1}, \dots, \X_k)
    g(\X_{k+1}, \dots, \X_{k+l-1}, \X, \X_{k+l+1}, \dots, \X_{2k}) \nonumber \\
    &\big| \Z = \z_j \; ; \; \Z_i=\z_i, \forall i=1, \dots, 2k, i \notin \{ j, k+l \} \big].
    \label{condUstat:eq:def:tilde_theta_jl} \\
    \theta_{j,l,m}&(\z_1, \dots, \z_{3k}) :=
    \EE \big[
    g(\X_1, \dots, \X_{j-1}, \X, \X_{j+1}, \dots, \X_k) \nonumber \\
    &g(\X_{k+1}, \dots, \X_{k+l-1}, \X, \X_{k+l+1}, \dots, \X_{2k})
    g(\X_{2k+1}, \dots, \X_{2k+m-1}, \X, \X_{2k+m+1}, \dots, \X_{3k}) \nonumber \\
    &\big| \Z = \z_j \; ; \; \Z_i=\z_i, \forall i=1, \dots, 3k, i \notin \{ j, k+l, 2k+m \} \big]. \nonumber
\end{align}

\begin{assumpt}
    \begin{enumerate}[(i)]
        \item $h_n \to 0$ and $n h_n^{p} \to \infty$ ;
        \item $K$ is symmetric at $0$, bounded and compactly supported ;
        \item $\theta_{j,l}$ is continuous at $(\z_1, \dots, \z_k)$ for all $1 \leq j,l \leq k$ ;
        \item $\theta$ is two times continuously differentiable in a neighborhood of $(\z_1, \dots, \z_k)$ ;
        \item $\theta_{j,l,m}$ is bounded in a neighborhood of $(\z_1, \dots, \z_k, \z_1, \dots, \z_k, \z_1, \dots, \z_k) \in \Zc^{3k}$, for all $1 \leq j,l,m \leq k$ ;
        \item $f_\Z$ is twice differentiable in neighborhoods of $\z_i, 1 \leq i \leq k.$
    \end{enumerate}
    \label{condUstat:assumpt:stute_asymptotic_normality}
\end{assumpt}

\begin{prop}[Asymptotic normality of $\hat \theta$, Corollary 2.4 in Stute \cite{stute1991conditional}]
    Under Assumption \ref{condUstat:assumpt:stute_asymptotic_normality},
    we have $\sqrt{n h_n^{p}} \big(\hat \theta(\z_1, \dots, \z_k) - \theta(\z_1, \dots, \z_k) \big) \indistrto \Nc(0,\rho^2),$ \\
    where $\rho^2 := \sum_{j,l = 1}^k \1_{\{\z_j = \z_l\}}
    \big(\theta_{j,l}(\z_1, \dots, \z_k) - \theta^2(\z_1, \dots, \z_k) \big)
    ||K||_2^2 / f_\Z(\z_j)$.
    
    \mds
    
    Moreover, let $N$ be a positive integer, and
    $\big(\z_1^{(1)}, \dots, \z_k^{(1)}, \dots,
    \z_1^{(N)}, \dots, \z_k^{(N)} \big) \in \Zc^{k \times N}$.
    Then under similar regularity conditions,
    $\sqrt{n h_n^{p}} \big( \hat \theta(\z_1^{(i)}, \dots, \z_k^{(i)})
    - \theta(\z_1^{(i)}, \dots, \z_k^{(i)})\big)_{i=1, \dots, N}
    \indistrto \Nc(0,\HH)$,
    where, for $1 \leq \tilde j, \tilde l \leq N$,
    \begin{equation*}
        [\HH]_{\tilde j, \tilde l} := \sum_{j,l = 1}^k
        \1_{\left\{\z_j^{(\tilde j)} = \z_l^{(\tilde l)}\right\}}
        \bigg(\tilde \theta_{j,l}
        \left(\z_1^{(\tilde j)}, \dots, \z_k^{(\tilde j)} ,
        \z_1^{(\tilde l)}, \dots, \z_k^{(\tilde l)} \right)
        - \theta \left(\z_1^{(\tilde j)}, \dots, \z_k^{(\tilde j)} \right) 
        \theta \left(\z_1^{(\tilde l)}, \dots, \z_k^{(\tilde l)} \right) \bigg)
        \frac{||K||_2^2}{f_\Z \left(\z_j^{(\tilde j)} \right)}.
    \end{equation*}
    \label{condUstat:prop:asympt_norm:hat_theta}
\end{prop}
Note that the second part of Proposition \ref{condUstat:prop:asympt_norm:hat_theta} above is a consequence of the first one. Indeed, for every $(c_1, \dots, c_{N}) \in \Rb^{N}$, we can define
$\theta \big(\z_1^{(1)}, \dots, \z_k^{(1)}, \dots,
\z_1^{(N)}, \dots, \z_k^{(N)} \big)
:= \sum_{\tilde i = 1}^{N} c_{\tilde i} \theta(\z_1^{(\tilde i)}, \dots, \z_k^{(\tilde i)})$ and corresponding versions of $g$, $\hat \theta$ and $\rho^2$.
Finally, the conclusion follows from the Cramér-Wold device.

\section{Finite distance proofs for \texorpdfstring{$\hat \theta$ and $\hat \beta$}{hat theta and hat beta}}

For convenience, we recall Berk's (1970) inequality (see Theorem A in Serfling \cite[p.201]{serfling1980approximation}). Note that, if $m=1$, this reduces to Bernstein's inequality.
\begin{lemma}
    Let $k > 0$, $n \geq k$, $\X_1, \dots, \X_n$ i.i.d. random vectors with values in a measurable space $\Xc$ and $g: \Xc^k \to [a, b]$ be a real bounded function. Set $\theta := \EE[g(\X_{1:k})]$ and $\sigma^2 := Var[g(\X_{1:k})]$.
    Then, for any $t > 0$,
    \begin{align*}
        \PP \left( \binom{n}{k}^{-1} \sum_{\sigma \in \Iknup}
        g \left(\X_{\sigma(1:k)} \right) - \theta \geq t \right)
        \leq \exp \bigg( - \frac{[ n / k] t^2}{2 \sigma^2 + (2/3) (b-\theta) t} \bigg),
    \end{align*}
    where $\Ikn$ is the set of bijective functions from $\{1, \dots, k\}$ to $\{1, \dots, n\}$ and $\Iknup$ is the subset of $\Ikn$ made of increasing functions.
    \label{lemma:bernstein_U_stat_sym}
\end{lemma}
Note that $g$ does not need to be symmetric for this bound to hold. Indeed, if $g$ is not symmetric, we can nonetheless apply this lemma to the symmetrized version $\tilde g$ defined as $\tilde g(\x_{1:k}) := (k!)^{-1} \sum_{\sigma \in \Ifrak{k}{k}} g(\x_{\sigma(1:k)}),$ and we get the result.

\subsection{Proof of Lemma \ref{lemma:bound_Nk_fk}}
\label{condUstat:proof:lemma:bound_Nk_fk}

We decompose the quantity to bound into a stochastic part and a bias as follows:
\begin{align*}
    N_k(\z_{1:k}) - \prod_{i=1}^k f_{\Z}(\z_i)
    = \big( N_k(\z_{1:k}) - \EE[N_k(\z_{1:k})] \big) 
    + \big( \EE[N_k(\z_{1:k})] - \prod_{i=1}^k f_{\Z}(\z_i) \big).
\end{align*}
We first bound the bias.
\begin{align*}
    \bigg| \EE \big[N_k(\z_{1:k}) \big] - \prod_{i=1}^k f_{\Z}(\z_i) \bigg|
    &= \bigg| \EE \Big[ \binom{n}{k}^{-1} \sum_{\sigma \in \Ikn} \prod_{i=1}^k
    K_h \big(\Z_{\sigma(i)}-\z_i \big) \Big] - \prod_{i=1}^k f_{\Z}(\z_i) \bigg| \\
    %
    &= \bigg| \int \Big( \prod_{i=1}^k f_{\Z}(\z_i + h \u_i)
    - \prod_{i=1}^k f_{\Z}(\z_i) \Big)
    \prod_{i=1}^k K ( \u_i ) d\u_i \bigg| \\
    &= \bigg| \int \Big(\phi_{\z, \u}(1) - \phi_{\z, \u}(0) \Big) 
    \prod_{i=1}^k K ( \u_i ) d\u_i \bigg|,
\end{align*}
where $\phi_{\z, \u}(t) := \prod_{j=1}^k f_{\Z} \big( \z_i + t h \u_j \big)$ for $t \in [-1,1]$. Note that this function has at least the same regularity as $f_\Z$, so it is $\alpha$-differentiable, and by a Taylor-Lagrange expansion, we get
\begin{align*}
    \bigg| &\EE[N_k(\z_{1:k})] - \prod_{i=1}^k f_{\Z}(\z_i) \bigg|
    = \bigg| \int_{\Rb^{k p}}
    \bigg(\sum_{i=1}^{\alpha-1} \frac{1}{i!} \phi_{\z, \u}^{(i)}(0) 
    + \frac{1}{\alpha !}\phi_{\z, \u}^{(\alpha)}
    (t_{\z, \u}) \bigg) \prod_{i=1}^k K ( \u_i ) d\u_i \bigg|.
\end{align*}
For $l>0$, we have
\begin{align*}
    \phi_{\z, \u}^{(l)}(0) 
    &= \sum_{m_1 + \, \cdots \, + m_k = l} \binom{\alpha}{m_{1:k}}
    \prod_{i=1}^k \frac{\partial^{m_i} \Big( f_{\Z} \big(\z_i + ht \u_i \big) \Big)}{\partial t^{m_i}}(0) \\
    &= \sum_{m_1 + \, \cdots \, + m_k = l} \binom{\alpha}{m_{1:k}}
    \prod_{i=1}^k \sum_{j_1, \dots, j_{m_i} = 1}^{p} h^{m_i} u_{i,j_1} \dots u_{i,j_{m_i}} 
    \frac{\partial^{m_i} f_{\Z}}{\partial z_{j_1} \, \cdots \, \partial z_{j_{m_i}}}
    \big(\z_i + t_{\z, \u} h \u_i \big),
\end{align*}
where $\binom{\alpha}{m_{1:k}} := \alpha! / \big( \prod_{i=1}^k (m_i!) \big)$ is the multinomial coefficient. Using Assumption \ref{condUstat:assumpt:kernel_integral}, for every
$i=1, \dots, \alpha - 1$, we get $\int K(\u_1) \cdots K(\u_k) \phi_{\z, \u}^{(i)}(0) d\u_1 \cdots d\u_k = 0$. Therefore, only the last term remains and we have
\begin{align*}
    \bigg| &\EE[N_k(\z_{1:k})] - \prod_{i=1}^k f_{\Z}(\z_i) \bigg|
    = \bigg| \int \bigg( \frac{1}{\alpha !}\phi_{\z, \u}^{(\alpha)}
    (t_{\z, \u}) \bigg) \prod_{i=1}^k K ( \u_i ) d\u_i \bigg|
    %
    \leq \frac{C_{K, \alpha}}{\alpha!} h^{\alpha},
\end{align*}
using Assumption \ref{condUstat:assumpt:bound_int_deriv_fZ}.

\medskip

Second, we bound the stochastic part.
We have
\begin{align*}
    N_k(\z_{1:k}) - \EE[N_k(\z_{1:k})]
    &= \frac{k!(n-k)!}{n!} \sum_{\sigma \in \Iknup}
    \prod_{i=1}^k K_h \big(\Z_{\sigma(i)}-\z_i \big)
    - \prod_{i=1}^k \EE[K_h \big(\Z_{i}-\z_i \big)].
\end{align*}
Then, we can apply Lemma \ref{lemma:bernstein_U_stat_sym} to the function $g$ defined by
$g(\tilde \z_1, \dots, \tilde \z_k) := \prod_{i=1}^k K_h \big(\tilde \z_i - \z_i \big)$.
Here, we have $b = -a = h^{-k p} C_K^k$, and
\begin{align*}
    Var \big[ g(\Z_1, \dots, \Z_k)^2 \big]
    \leq \EE \big[ g(\Z_1, \dots, \Z_k)^2 \big]
    = \prod_{i=1}^k \EE \big[ K_h \big( \Z_i - \z_i \big)^2 \big]
    \leq h^{-k p} f_{\Z, max}^k ||K||_2^{2k}.
\end{align*}
Finally, we get
\begin{align*}
    \PP \left( \binom{n}{k}^{-1} N_k(\z_{1:k}) - \EE[N_k(\z_{1:k})] \geq t \right)
    \leq \exp \bigg( - \frac{[ n / k] t^2}{2 h^{-k p} f_{\Z, max}^k ||K||_2^{2k} + (4/3) h^{-k p} C_K^k t} \bigg),
\end{align*}

\begin{flushright}
    $\Box$
\end{flushright}

\subsection{Proof of Proposition \ref{condUstat:prop:exponential_bound_hat_theta}}
\label{condUstat:proof:prop:exponential_bound_hat_theta}

We have the following decomposition
\begingroup \allowdisplaybreaks
\begin{align*}
    | &\hat \theta(\z_{1:k}) - \theta(\z_{1:k}) | \\
    &= \bigg| N_k(\z_{1:k})^{-1} \frac{(n-k)!}{n!}
    \sum_{\sigma \in \Ikn} \prod_{i=1}^k K_h \big(\Z_{\sigma(i)}-\z_i \big)
    \Big( g(\X_{\sigma(1:k)}) - \EE \big[ g(\X_{1:k}) \big| \Z_{1:k}=\z_{1:k} \big] \Big) \bigg| \\
    &= \frac{\prod_{i=1}^k f_{\Z}(\z_i)}{N_k(\z_1, \dots, \z_k)} \cdot
    \bigg| \frac{(n-k)!}{n!} \sum_{\sigma \in \Ikn}
    \prod_{i=1}^k \frac{K_h \big(\Z_{\sigma(i)}-\z_i \big)}{f_{\Z}(\z_i)}
    \Big( g(\X_{\sigma(1:k)}) - \EE \big[ g(\X_{1:k}) \big| \Z_{1:k}=\z_{1:k} \big] \Big) \bigg| \\
    &=: \frac{\prod_{i=1}^k f_{\Z}(\z_i)}{N_k(\z_1, \dots, \z_k)} \; \cdot
    \bigg| \sum_{\sigma \in \Ikn} S_{\sigma} \bigg|.
\end{align*}
\endgroup
The conclusion will follow from the next three lemmas, where we will bound separately $\prod_{i=1}^k f_{\Z} / N_k$,
the bias term $\big| \sum_{\sigma \in \Ikn} \EE[S_{\sigma}] \big|$
and the stochastic component $\big| \sum_{\sigma \in \Ikn} \big(S_{\sigma} - \EE[S_{\sigma}] \big) \big|$.

\begin{lemma}[Bound for $\prod_{i=1}^k f_{\Z}(\z_i)/N_k$]
    Under Assumptions \ref{condUstat:assumpt:kernel_integral},
    \ref{condUstat:assumpt:bound_int_deriv_fZ}, \ref{condUstat:assumpt:f_Z_max}, and \ref{condUstat:assumpt:f_Z_min}
    and if for some $t > 0$, $C_{K, \alpha} h^{\alpha} / \alpha ! + t < f_{\Z, min}^k/2$,
    we have
    \begin{align*}
        &\PP \Bigg( \bigg| \frac{1}{N_k(\z_{1:k})} - \frac{1}{\prod_{i=1}^k f_{\Z}(\z_i)} \bigg|
        \leq \frac{4}{f_{\Z, min}^{2k}} \bigg( \frac{C_{K, \alpha}  h^{\alpha}}{\alpha !} + t \bigg)
        \Bigg) \\
        &\hspace{6cm} \geq 1 - 2 \exp \bigg( - \frac{[ n / k] t^2}{2 h^{-k p} f_{\Z, max}^k ||K||_2^{2k} + (4/3) h^{-k p} C_K^k t} \bigg),
    \end{align*}
    and on the same event, $N_k(\z_{1:k})$ is strictly positive and
    \begin{align*}
        \frac{\prod_{i=1}^k f_{\Z}(\z_i)}{N_k(\z_{1:k})}
        \leq 1 + \frac{4 f_{\Z, max}^k}{f_{\Z, min}^{2k}} \bigg( \frac{C_{K, \alpha} h^{\alpha}}{\alpha !} + t \bigg).
    \end{align*}
    \label{lemma:bound_prodf_Nk}
\end{lemma}

{\it Proof :} Using the mean value inequality for the function $x \mapsto 1/x$, we get
\begin{align*}
    \bigg| \frac{1}{N_k(\z_{1:k})}-\frac{1}{\prod_{i=1}^k f_{\Z}(\z_i)} \bigg|
    \leq \frac{1}{N_*^2} \big|N_k(\z_{1:k})-\prod_{i=1}^k f_{\Z}(\z_i) \big|,
\end{align*}
where $N_*$ lies between $N_k(\z_{1:k})$ and $\prod_{i=1}^k f_{\Z}(\z_i)$.
By Lemma \ref{lemma:bound_Nk_fk}, we get 
\begin{align*}
        \PP \bigg( \big| N_k(\z_{1:k}) &- \prod_{i=1}^k f_{\Z}(\z_i) \big| \leq \frac{C_{K, \alpha}}{\alpha !} h^{\alpha} + t \bigg)
        \geq 1 - 2 \exp \bigg( - \frac{[ n / k] t^2}{2 h^{-k p} f_{\Z, max}^k ||K||_2^{2k} + (4/3) h^{-k p} C_K^k t} \bigg).
    \end{align*}
On this event,
$\big| N_k(\z_{1:k}) - \prod_{i=1}^k f_{\Z}(\z_i) \big| \leq (1/2) \prod_{i=1}^k f_{\Z}(\z_i)$ by assumption,
so that $f_{\Z, min}^k/2 \leq N_k(\z_{1:k})$. We have also $f_{\Z, min}^k/2 \leq \prod_{i=1}^k f_{\Z}(\z_i)$. Thus, we have $f_{\Z, min}^k/2 \leq N_*$.
Combining the previous inequalities, we finally get
\begin{align*}
    \bigg| \frac{1}{N_k(\z_{1:k})}
    - \frac{1}{\prod_{i=1}^k f_{\Z}(\z_i)} \bigg|
    \leq \frac{1}{N_*^2} \big| N_k(\z_{1:k})
    - \prod_{i=1}^k f_{\Z}(\z_i) \big|
    \leq \frac{4}{f_{\Z, min}^{2k}} \bigg( \frac{C_{K, \alpha} h^{\alpha}}{\alpha !} + t \bigg).
\end{align*}

\begin{flushright}
    $\Box$
\end{flushright}

Now, we provide a bound on the bias.

\begin{lemma}
    Under Assumptions \ref{condUstat:assumpt:kernel_integral} and \ref{condUstat:assumpt:f_XZ_Holder}, we have
    $\big| \EE[S_{\sigma}] \big|
    \leq C_{g, f, \alpha} C_{K,\alpha} h^{k \alpha} / (f_{\Z,min}^k \alpha !).$
    \label{lemma:boundProba_sum_S_sigma:bias}
\end{lemma}
\noindent
{\it Proof :} We remark that
\begin{align}
    0 &= \int \Big( g(\x_{1:k})
    - \EE \big[ g(\X_{1:k}) \big| \Z_{1:k} = \z_{1:k} \big] \Big)
    f_{\X | \Z=\z_1}(\x_1) \cdots f_{\X | \Z=\z_k}(\x_k) d\mu^{\otimes k}(\x_{1:k}) \nonumber \\
    &= \int \Big( g(\x_{1:k})
    - \EE \big[ g(\X_{1:k}) \big| \Z_{1:k} = \z_{1:k} \big] \Big)
    \frac{f_{\X,\Z}(\x_1,\z_1)
    \cdots f_{\X, \Z}(\x_k,\z_k)}{\prod_{i=1}^k f_{\Z}(\z_i)} d\mu^{\otimes k}(\x_{1:k}).
\end{align}
We have
\begin{align*}
    &\EE[S_{\sigma}]
    = \EE \bigg[\frac{K_h(\Z_{\sigma(1)}-\z_1) \cdots K_h(\Z_{\sigma(k)}-\z_k)}{\prod_{i=1}^k f_{\Z}(\z_i)}
    \Big( g\big(\X_{\sigma(1)}, \dots, \X_{\sigma(k)} \big)
    - \EE \big[ g(\X_{1:k}) \big| \Z_{1:k} = \z_{1:k} \big] \Big) \bigg] \displaybreak[0] \\
    %
    %
    &= \int \Big( g(\x_{1:k})
    - \EE \big[ g(\X_{1:k}) \big| \Z_{1:k} = \z_{1:k} \big] \Big) 
    \prod_{i=1}^k \frac{K(\u_i)}{f_{\Z}(\z_i)} f_{\X,\Z}(\x_i, \z_i + h\u_i) \,
    d\mu(\x_i) d\u_i \displaybreak[0] \\
    &= \int\Big( g(\x_{1:k}) - \EE \big[ g(\X_{1:k}) \big| \Z_{1:k} = \z_{1:k} \big] \Big) 
    \bigg( \prod_{i=1}^k f_{\X,\Z} \Big(\x_i, \z_i + h \u_i \Big) 
    - \prod_{i=1}^k f_{\X,\Z} \Big(\x_i, \z_i \Big) \bigg) \,
    \prod_{i=1}^k \frac{K(\u_i)}{f_{\Z}(\z_i)} \, d\mu(\x_i) d\u_i.
\end{align*}
We apply now the Taylor-Lagrange formula to the function
\begin{equation*}
    \phi_{\x_{1:k}, \u_{1:k}}(t)
    := \prod_{i=1}^k f_{\X,\Z} \Big(\x_i, \z_i + h \u_i \Big),
\end{equation*}
and get
\begin{align*}
    &\EE[S_{\sigma}]
    = \int \Big( g(\x_{1:k})
    - \EE \big[ g(\X_{1:k}) \big| \Z_{1:k} = \z_{1:k} \big] \Big)
    \bigg( \phi_{\x_{1:k}, \u_{1:k}}(t)(1)
    - \phi_{\x_{1:k}, \u_{1:k}}(t)(0) \bigg) 
    \prod_{i=1}^k \frac{K(\u_i)}{f_{\Z}(\z_i)} \, d\mu(\x_i) d\u_i \\
    &= \int \Big( g(\x_{1:k})
    - \EE \big[ g(\X_{1:k}) \big| \Z_{1:k} = \z_{1:k} \big] \Big)  \\
    & \hspace{3cm} \cdot
    \bigg( \sum_{j=1}^{\alpha - 1} \frac{1}{j !}
    \phi_{\x_{1:k}, \u_{1:k}}(t)^{(j)}(0) + \frac{1}{\alpha !}
    \phi_{\x_{1:k}, \u_{1:k}}(t)^{(\alpha)}
    (t_{\x, \u}) \bigg) \,
    \prod_{i=1}^k \frac{K(\u_i)}{f_{\Z}(\z_i)} \, d\mu(\x_i) d\u_i 
    \displaybreak[0] \\
    &= \int \Big( g(\x_{1:k})
    - \EE \big[ g(\X_{1:k}) \big| \Z_{1:k} = \z_{1:k} \big] \Big) \\
    & \hspace{3cm} \cdot
    \bigg( \frac{1}{\alpha !}
    \phi_{\x_{1:k}, \u_{1:k}}(t)^{(\alpha)} (t_{\x, \u}) \bigg) \,
    \prod_{i=1}^k \frac{K(\u_i)}{f_{\Z}(\z_i)} \, d\mu(\x_i) d\u_i
    \displaybreak[0] \\
    &= \int \Big( g(\x_{1:k})
    - \EE \big[ g(\X_{1:k}) \big| \Z_{1:k} = \z_{1:k} \big] \Big) \\
    & \hspace{3cm} \cdot
    \frac{1}{\alpha !} \bigg(
    \phi_{\x_{1:k}, \u_{1:k}}(t)^{(\alpha)} (t_{\x, \u})
    - \phi_{\x_{1:k}, \u_{1:k}}(t)^{(\alpha)} (0) \bigg) \,
    \prod_{i=1}^k \frac{K(\u_i)}{f_{\Z}(\z_i)} \, d\mu(\x_i) d\u_i.
\end{align*}
For every real $t$, we have
\begin{align}
    &\phi^{(\alpha)} (t)
    = \sum_{m_1 + \, \cdots \, + m_k = \alpha} \binom{n}{m_{1:k}}
    \prod_{i=1}^k \frac{\partial^{m_i} \Big( f_{\X,\Z} \big(\x_i,\z_i + ht \u_i \big) \Big)}{\partial t^{m_i}} 
    \nonumber \\
    &= \sum_{m_1 + \, \cdots \, + m_k = \alpha} \binom{n}{m_{1:k}}
    \prod_{i=1}^k \sum_{j_1, \dots, j_{m_i} = 1}^{p} h^{m_i} u_{i,j_1} \dots u_{i,j_{m_i}} 
    \frac{\partial^{m_i} f_{\X,\Z}}{\partial z_{j_1} \, \cdots \, \partial z_{j_{m_i}}} \big(\x_i,\z_i + ht \u_i \big)
    \nonumber \\
    &= h^{\alpha}
    \sum_{m_1 + \, \cdots \, + m_k = \alpha} \binom{n}{m_{1:k}}
    \prod_{i=1}^k \sum_{j_1, \dots, j_{m_i} = 1}^{p} u_{i,j_1} \dots u_{i,j_{m_i}} 
    \frac{\partial^{m_i} f_{\X,\Z}}{\partial z_{j_1} \, \cdots \, \partial z_{j_{m_i}}} \big(\x_i,\z_i + ht \u_i \big).
    \label{condUstat:eq:computation_partial_deriv_phi}
\end{align}
Therefore, we get
\begin{align*}
    &\EE[S_{\sigma}]
    = \sum_{m_1 + \, \cdots \, + m_k = \alpha}
    \binom{n}{m_{1:k}}
    \int \prod_{i=1}^k \frac{K(\u_i)}{\prod_{i=1}^k f_{\Z}(\z_i)}
    \sum_{j_1, \dots, j_{m_i} = 1}^{p} u_{i,j_1} \dots u_{i,j_{m_i}} \\
    &\hspace{3cm} \cdot \Big( g(\x_{1:k})
    - \EE \big[ g(\X_{1:k}) \big| \Z_{1:k} = \z_{1:k} \big] \Big) \\
    &\hspace{3cm} \cdot \left( \frac{\partial^{m_i} f_{\X,\Z}}{\partial z_{j_1} \, \cdots \, \partial z_{j_{m_i}}} \big(\x_i,\z_i + ht \u_i \big) - \frac{\partial^{m_i} f_{\X,\Z}}{\partial z_{j_1} \, \cdots \, \partial z_{j_{m_i}}} \big(\x_i,\z_i \big) \right) d\mu(\x_1) d\u_1 \cdots d\mu(\x_k) d\u_k,
\end{align*}
and, using Assumption~\ref{condUstat:assumpt:f_XZ_Holder}, this yields
\begin{align*}
    \big| \EE[S_{\sigma}] \big|
    & \leq \frac{C_{g, f, \alpha} C_{K,\alpha} h^{\alpha+k}}
    {f_{\Z,min}^k \alpha !}.
\end{align*}

\begin{flushright}
    $\Box$
\end{flushright}

\mds

Now we bound the stochastic component. We have the following equality
\begin{align*}
    \bigg| \sum_{\sigma \in \Ikn} \big( S_{\sigma} - \EE[S_{\sigma}] \big) \bigg|
    = \bigg| \frac{(n-k) !}{n !}
    \sum_{\sigma \in \Ikn} g\big( (\X_{\sigma(1)}, \Z_{\sigma(1)}) \, , \dots , \, 
    (\X_{\sigma(k)}, \Z_{\sigma(k)}) \big) \bigg|
\end{align*}
with the function $\tilde g$ defined by
\begin{align*}
    &\tilde g\big( (\X_1, \Z_1) \, , \dots , \, (\X_k, \Z_k) \big) \\
    &= \frac{K_h \big(\Z_1-\z_1 \big) \cdots K_h \big(\Z_k-\z_k \big)}{\prod_{i=1}^k f_{\Z}(\z_i)}
    \Big( g(\X_{1:k}) - \EE \big[ g(\X_{1:k}) \big| \Z_{1:k} = \z_{1:k} \big] \Big) \\
    & - \EE \Bigg[ \frac{K_h \big(\Z_1-\z_1 \big) \cdots K_h \big(\Z_k-\z_k \big)}{\prod_{i=1}^k f_{\Z}(\z_i)}
    \Big( g(\X_{1:k}) - \EE \big[ g(\X_{1:k}) \big| \Z_{1:k} = \z_{1:k} \big] \Big) \Bigg]
\end{align*}
By construction, $\EE \Big[ \tilde g\big( (\X_1, \Z_1) \, , \dots , \, (\X_k, \Z_k) \big) \Big] = 0.$
If $\tilde g$ is bounded, we can derive an immediate bound for this stochastic component.
Indeed, we would have $||\tilde g||_\infty \leq 4 C_K^k h^{-kp} C_{g}^k / f_{\Z,min}^k$. 
Moreover, we have
\begin{align*}
    Var \Big[ \tilde g \big( (\X_1, \Z_1) \, , \dots , \, (\X_k, \Z_k) \big) \Big]
    &\leq \EE \bigg[ \frac{K_h^2 \big(\Z_1-\z_1 \big) \cdots K_h^2 \big(\Z_k-\z_k \big)}{\prod_{i=1}^k f_{\Z}^2(\z_i)}
    g{}^2(\X_1, \dots, \X_k) \bigg] \\
    &\leq C_{g}^2 f_{\Z,max}^k f_{\Z,min}^{-2k} h^{-k p} ||K||_2^{2k}.
\end{align*}
Therefore, we can apply Lemma \ref{lemma:bernstein_U_stat_sym}, and we get
\begin{align*}
    \PP \bigg( \Big| \sum_{\sigma \in \Ikn} &\big( S_{\sigma} - \EE[S_{\sigma}] \big) \Big| > t \bigg) 
    \leq 2 \exp \bigg( - \frac{[ n / k] t^2}
    {2 C_{g}^2 f_{\Z,max}^k f_{\Z,min}^{-2k} h^{-k p} ||K||_2^{2k} + (8/3) C_K^k h^{-kp} C_{g}^k f_{\Z,min}^{-k} t} \bigg).
\end{align*}

\mds 

In the following Lemma \ref{lemma:boundProba_sum_S_sigma}, our goal will be to bound the stochastic component using only Assumption~\ref{condUstat:assumpt:conditional_Bernstein} on the conditional moments of $g$.

\begin{lemma}
    Under Assumptions \ref{condUstat:assumpt:kernel_integral}, \ref{condUstat:assumpt:f_Z_min} and \ref{condUstat:assumpt:conditional_Bernstein}, for every $t > 0$, we have
    \begin{align*}
        \PP \left( \sum_{\sigma \in \Ikn} S_{\sigma} - \EE[S_{\sigma}] > t \right)
        &\leq \exp \left( - \frac{t^2 f_{\Z, min}^{2k} h^{k p} [ n / k]}{128 \big(B_{g, \z} + \tilde B_{g} \big)^2 C_K^{2k-1}
        + 2 t \big(B_{g, \z} + \tilde B_{g} \big) C_K^{k} f_{\Z, min}^{k}} \right).
    \end{align*}
    \label{lemma:boundProba_sum_S_sigma}
\end{lemma}

\noindent
{\it Proof:}
Using the same decomposition for U-statistics as in Hoeffding \cite{hoeffding1963probability}, we obtain
\begin{align*}
    \sum_{\sigma \in \Ikn} S_{\sigma} - \EE[S_{\sigma}]
    = \frac{1}{n!} \sum_{\sigma \in \Inn} \frac{1}{[ n / k]}
    \sum_{i = 1}^{[ n / k]} V_{n,i,\sigma},
\end{align*}
where
\begin{align*}
    V_{n,i,\sigma}
    &:= \tilde g \Big( \big(\X_{\sigma(1 + (i-1) k)}, \Z_{\sigma(2 + (i-1) k)} \big), \dots, \big(\X_{\sigma(i k)}, \Z_{\sigma(j k)} \big) \Big).
\end{align*}

\mds

For any $\lambda > 0$, we have
\begingroup \allowdisplaybreaks
\begin{align}
    \PP \left( \sum_{\sigma \in \Ikn} S_{\sigma} - \EE[S_{\sigma}] > t \right)
    &\leq e^{-\lambda t} \EE \left[ \exp \left( \lambda \sum_{\sigma \in \Ikn} S_{\sigma} - \EE[S_{\sigma}] \right) \right] \nonumber \\
    &\leq e^{-\lambda t} \EE \left[ \exp \left( \lambda
    \frac{1}{n!} \sum_{\sigma \in \Inn} \frac{1}{[ n / k]}
    \sum_{i = 1}^{[ n / k]} V_{n,i,\sigma} \right) \right] \nonumber \\
    &\leq  e^{-\lambda t} \frac{1}{n!} \sum_{\sigma \in \Inn} 
    \EE \left[ \exp \left( \lambda \frac{1}{[ n / k]}
    \sum_{i = 1}^{[ n / k]} V_{n,i,\sigma} \right) \right] \nonumber \\
    &\leq  e^{-\lambda t} \frac{1}{n!} \sum_{\sigma \in \Inn}
    \prod_{i = 1}^{[ n / k]}
    \EE \left[ \exp \left( \lambda \frac{1}{[ n / k]}
    V_{n,i,\sigma} \right) \right] \nonumber \\
    &\leq  e^{-\lambda t} \left(
    \sup_{\sigma \in \Inn, \, i = 1, \dots, [ n / k]}
    \EE \left[ \exp \left( \lambda [ n / k]^{-1}
    V_{n,i,\sigma} \right) \right]
    \right)^{[ n / k]}.
    \label{condUstat:eq:proof:beggining_Hoeffding}
\end{align}
\endgroup

Let $l \geq 2$. Using the inequality
$(a+b+c+d)^l \leq 4^l (a^l+b^l+c^l+d^l)$, we get
\begin{align*}
    \EE \big[ |V_{n,i,\sigma}|^l \big]
    = \EE \big[ |V_{n,1,\sigma}|^l \big]
    & \leq 4^l \,
    \EE \Big[ |g(\X_{\sigma(1)}, \dots, \X_{\sigma(k)}) |^l
    \prod_{i=1}^{k} \frac{|K_h|^l \big(\Z_{\sigma(i)}-\z_i \big)}
    {f_{\Z}^l(\z_i)} \Big] \\
    & \hspace{0cm} + 4^l \,
    \EE \Big[ \big| \EE \big[ g(\X_{1:k})
    \big| \Z_{1:k} = \z_{1:k} \big] \big|^l
    \prod_{i=1}^{k} \frac{|K_h|^l \big(\Z_{\sigma(i)}-\z_i \big)}{f_{\Z}^l(\z_i)} \Big] \\
    & + 4^l \,
    \bigg| \EE \Big[g(\X_{\sigma(1)}, \dots, \X_{\sigma(k)})
    \prod_{i=1}^{k} \frac{K_h \big(\Z_{\sigma(i)}-\z_i \big)}
    {f_{\Z}^l(\z_i)} \Big] \bigg|^l \\
    & \hspace{0cm} + 4^l \,
    \bigg| \EE \Big[ \big| \EE \big[ g(\X_{1:k}) \big| 
    \Z_{1:k} = \z_{1:k} \big] \big| 
    \prod_{i=1}^{k} \frac{K_h \big(\Z_{\sigma(i)}-\z_i \big)}{f_{\Z}^l(\z_i)} \Big] \bigg|^l
\end{align*}
Using Jensen's inequality for the function $x \mapsto |x|^p$ with the second, third and fourth terms, and the law of iterated expectations for the first and the third terms, we get
\begin{align}
    \EE \big[ |V_{n,i,\sigma}|^l \big]
    & \leq 4^l \cdot 2 \,\EE \Big[ \EE \big[ |g(\X_{\sigma(1)}, \dots, \X_{\sigma(k)}) |^l \big| \Z_{\sigma(1)}, \dots, \Z_{\sigma(k)} \big] \prod_{i=1}^{k} \frac{|K_h|^l \big(\Z_{\sigma(i)}-\z_i \big)}
    {f_{\Z}^l(\z_i)} \Big] \nonumber \\
    & \hspace{0cm} + 4^l \cdot 2 \,\EE \Big[
    \EE \big[ \big| g(\X_{1:k}) \big|^l
    \big| \Z_i=\z_i, \forall i=1, \dots, k \big] \prod_{i=1}^{k} \frac{|K_h|^l \big(\Z_{\sigma(i)}-\z_i \big)}{f_{\Z}^l(\z_i)} \Big] \displaybreak[0] \nonumber \\
    & \leq 4^l \cdot 2 \, \EE \Big[ \Big(
    B_{g}^l(\Z_1, \dots, \Z_k) + B_{g}^l(\z_{1}, \dots, \z_{k})
    \Big)^l l!
    \prod_{i=1}^{k} \frac{|K_h|^l \big(\Z_{\sigma(i)}-\z_i \big)}{f_{\Z}^l(\z_i)} \Big] \displaybreak[0] \nonumber \\
    & \leq 4^l \cdot 2 \, \big(
    \tilde B_{g}^l + B_{g}^l(\z_{1}, \dots, \z_{k})
    \big) l! (h^{-k p} C_K^{k} f_{\Z, min}^{-k})^{l-1} \, f_{\Z, min}^{-k}
    \displaybreak[0] \nonumber \\
    %
    & \leq 2 \Big(4 \big(\tilde B_{g} + B_{g, \z} \big)
    h^{-k p} C_K^{k} f_{\Z, min}^{-k} \Big)^l l! \, h^{k p} C_K^{-1} \, , \nonumber
\end{align}
where $B_{g, \z} := B_{g}(\z_1, \dots, \z_k)$.
Remarking that $\EE[V_{n,i,\sigma}]=0$ by construction of $\tilde g$, we obtain
\begingroup \allowdisplaybreaks
\begin{align*}
    \EE \left[ \exp \left( \lambda [ n / k]^{-1} V_{n,i,\sigma} \right) \right]
    &= 1 + \sum_{l=2}^{\infty}
    \frac{\EE \big[ (\lambda [ n / k]^{-1} V_{n,i,\sigma})^l \big]}{l!} \\
    &\leq 1 + 2  C_K^{-1} h^{k p} \sum_{l=2}^{\infty}
    (4 \lambda [ n / k]^{-1} \big(B_{g, \z} + \tilde B_{g} \big) h^{-k p} C_K^{k} f_{\Z, min}^{-k})^l \\
    &\leq 1 + 2  C_K^{-1} h^{k p} \cdot \frac{ \left(4 \lambda [ n / k]^{-1}
    \big(B_{g, \z} + \tilde B_{g} \big) h^{-k p} C_K^{k} f_{\Z, min}^{-k} \right)^2}
    {1 - 4 \lambda [ n / k]^{-1} \big(B_{g, \z} + \tilde B_{g} \big) h^{-k p} C_K^{k} f_{\Z, min}^{-k}} \\
    &\leq \exp \left( \frac{
    32 \lambda^2 [ n / k]^{-2} \big(B_{g, \z} + \tilde B_{g} \big)^2 h^{- k p} C_K^{2k-1} f_{\Z, min}^{-2k}}
    {1 - 4 \lambda [ n / k]^{-1} \big(B_{g, \z} + \tilde B_{g} \big) h^{-k p} C_K^{k} f_{\Z, min}^{-k}} \right),
\end{align*}
\endgroup
where the last statement follows from the inequality $1+x \leq \exp (x)$.
Combining the latter bound with Equation (\ref{condUstat:eq:proof:beggining_Hoeffding}), we get
\begin{align}
    \PP \left( \sum_{\sigma \in \Ikn} S_{\sigma} - \EE[S_{\sigma}] > t \right)
    & \leq
    \exp \left( -\lambda t + \frac{32 \lambda^2 \big(B_{g, \z} + \tilde B_{g} \big)^2 C_K^{2k-1}}
    {f_{\Z, min}^{2k} h^{k p} [ n / k] - 4 \lambda \big(B_{g, \z} + \tilde B_{g} \big) C_K^{k} f_{\Z, min}^{k}}
    \right).
    \label{condUstat:eq:proof:Hoeffding_before_minimization}
\end{align}
Remarking that the right-hand side term inside the exponential is of the form $-\lambda t + \frac{a \lambda^2}{b - c \lambda}$, we choose the value
\begin{align}
    \lambda_* = \frac{tb}{2a+tc}
    = \frac{t f_{\Z, min}^{2k} h^{k p} [ n / k]}
    {64 \big(B_{g, \z} + \tilde B_{g} \big)^2  C_K^{2k-1}
    + t \big(B_{g, \z} + \tilde B_{g} \big) C_K^{k} f_{\Z, min}^{k}}
    \label{condUstat:eq:def:lambda_star_Bernstein}
\end{align}
such that $-\lambda_* t + \frac{a \lambda_*^2}{b - c \lambda_*} = - \frac{t^2 b}{4a + 2ct} = - \frac{t}{2} \lambda_*$.
Therefore, the right-hand side term of Equation (\ref{condUstat:eq:proof:Hoeffding_before_minimization}) can be simplified, and combining this with Equation (\ref{condUstat:eq:def:lambda_star_Bernstein}), we obtain
\begin{align*}
    \PP \left( \sum_{\sigma \in \Ikn} S_{\sigma} - \EE[S_{\sigma}] > t \right)
    &\leq \exp \left( - \frac{t^2 f_{\Z, min}^{2k} h^{k p} [ n / k]}{128 \big(B_{g, \z} + \tilde B_{g} \big)^2 C_K^{2k-1}
    + 2 t \big(B_{g, \z} + \tilde B_{g} \big) C_K^{k} f_{\Z, min}^{k}} \right).
\end{align*}

\begin{flushright}
    $\Box$
\end{flushright}

\subsection{Proof of Theorem \ref{condUstat:thm:bound_proba_hat_beta} }
\label{condUstat:proof:thm:bound_proba_hat_beta}

By Proposition \ref{condUstat:prop:exponential_bound_hat_theta}, for every $t_1, t_2 > 0$ such that $C_{K, \alpha} h^{\alpha} / \alpha ! + t < f_{\Z, min}/2$, we have
\begin{align*}
    &\PP \Bigg( |\hat \theta(\z_1, \dots, \z_k) - \theta(\z_1, \dots, \z_k) |
    < \big(1 + C_3 h^{\alpha} + C_4 t_1 \big) \times
    \big( C_5 h^{k + \alpha} + t_2 \big) \Bigg) \\
    & \hspace{4cm} \geq 1
    - 2 \exp \bigg( - \frac{[ n / k] t_1^2 h^{k p}}{ C_1 + C_2 t_1} \bigg)
    - 2 \exp \bigg( - \frac{[ n / k] t_2^2 h^{kp}}
    {C_6 + C_7 t_2} \bigg),
\end{align*}

We apply this proposition to every $k$-tuple $\big(\z'_{\sigma(1)}, \dots, \z'_{\sigma(k)} \big)$ where $\sigma \in \Iknpr$. Combining it with Assumption \ref{condUstat:assumpt:compact_case}, we get
\begin{align*}
    &\PP \Bigg( \sup_i |\xi_{i,n}| < C_{\Lambda'}
    \big(1 + C_3 h^{\alpha} + C_4 t_1 \big) \times
    \big( C_5 h^{k + \alpha} + t_2 \big) \Bigg) \\
    & \hspace{4cm} \geq 1
    - 2 \sum_{i = 1}^{|\Iknpr|} \Bigg[ \exp \bigg( - \frac{[ n / k] t_1^2 h^{k p}}{ C_1 + C_2 t_1} \bigg)
    + \exp \bigg( - \frac{[ n / k] t_2^2 h^{kp}}
    {C_6 + C_7 t_2} \bigg) \Bigg] ,
\end{align*}
Choosing $t_1 := f_{\Z, min}/4$ and using the bound (\ref{condUstat:cond:h_1_2}) on $h$, we get
\begin{align*}
    &\PP \Bigg( \sup_i |\xi_{i,n}| < C_{\Lambda'}
    \big(1 + C_3 \frac{f_{\Z, min} \alpha!}{4 C_{K, \alpha}} + C_4 \frac{f_{\Z, min}}{4} \big) \times
    \big( C_5 h^{k + \alpha} + t_2 \big) \Bigg) \\
    & \hspace{2cm} \geq 1
    - 2 \sum_{i = 1}^{|\Iknpr|} \Bigg[
    \exp \bigg( - \frac{[ n / k] f_{\Z, min}^2 h^{k p}}
    { 16 C_1 + 4 C_2 f_{\Z, min} } \bigg)
    + \exp \bigg( - \frac{[ n / k] t_2^2 h^{kp}}
    {C_6 + C_7 t_2} \bigg) \Bigg] .
\end{align*}
Choosing $t_2 = t / (2 C_8) = t / \Big( 2 C_\psibm C_{\Lambda'}
\big(1 + C_3 \frac{f_{\Z, min} \alpha!}{4 C_{K, \alpha}} + C_4 \frac{f_{\Z, min}}{4} \big) \Big) $, and using the bound (\ref{condUstat:cond:h_1_2}) on $h^\alpha$, we get
\begin{align*}
    &\PP \Bigg( \sup_i |\xi_{i,n}| < t / C_{\psibm} \Bigg) \geq 1 - 2 \sum_{i = 1}^{|\Iknpr|} \Bigg[
    \exp \bigg( - \frac{[ n / k] f_{\Z, min}^2 h^{k p}}
    { 16 C_1 + 4 C_2 f_{\Z, min} } \bigg)
    + \exp \bigg( - \frac{[ n / k] t^2 h^{kp}}
    {4 C_8^2 C_6 + 2 C_8 C_7 t} \bigg) \Bigg].
\end{align*}

On the same event, we have $\max_{j=1, \dots, p'} \Big| \frac{1}{n'}
\sum_{i=1}^{n'} Z'_{i,j} \xi_{i,n} \Big| \leq t$, by Assumption \ref{condUstat:assumpt:compact_case}.
The conclusion results from the following lemma.

\begin{lemma}[From {\protect\cite[Lemma 25]{derumigny2018kendall}}]
    Assume that
    $\max_{j=1, \dots, p'} \Big| \frac{1}{n'}
    \sum_{i=1}^{n'} Z'_{i,j} \xi_{i,n} \Big| \leq t,$
    for some $t > 0$, that the assumption $RE(s,3)$ is satisfied, and that the tuning parameter is given by $\lambda = \gamma t$, with $\gamma \geq 4$.
    Then, $||\Zb'(\hat \beta - \beta^*)||
    \leq \dfrac{4(\gamma+1) t \sqrt{s}}{\kappa(s,3)}$
    and $|\hat \beta - \beta^*|_q
    \leq \dfrac{4^{2/q}(\gamma+1) t s^{1/q} }{\kappa^2(s,3)},$
    for every $1 \leq q \leq 2$.
    \label{lemma:lasso_bound}
\end{lemma}

\begin{flushright}
    $\Box$
\end{flushright}

\section{Proof of Theorem \ref{condUstat:thm:weak_conv_doubleAsympt}}
\label{condUstat:proof:thm:weak_conv_doubleAsympt}

We detail the assumption which we will use to prove Theorem \ref{condUstat:thm:weak_conv_doubleAsympt}.

\begin{assumpt}
\begin{enumerate}[(i)]
    \item \label{condUstat:assumpt:kernel_separation_Zprime}
    The support of the kernel $K(\cdot)$ is included into $[-1,1]^{p}$.
    Moreover, for all $n,n'$ and every $(i, j) \in \{1, \dots, n'\}^2$, $i \neq j$, we have
    $|\z'_{i} - \z'_{j} |_\infty > 2 h_{n,n'}$.
    \item
    (a) $n' (n h_{n,n'}^{p+4\alpha} + h_{n,n'}^{2\alpha} + h_{n,n'}^{p} + (n h_{n,n'}^{p})^{-1} ) \to 0$, (b) $\lambda_{n,n'}  ( n' \, n \,  h_{n,n'}^{p})^{1/2} \to 0$,
    (c) $ n' \, n \,  h_{n,n'}^{p} \to \infty$ and $ n \,  h_{n,n'}^{p + 2 \alpha - \epsilon}/\ln n' \to \infty$ for some $\epsilon \in [ 0, 2\alpha[$.
    \item \label{condUstat:assumpt:limit_distrib_Zprime}
    The distribution
    $\PP_{\z',n'} := |\Iknpr|^{-1} \sum_{\sigma \in \Iknpr} 
    \delta_{(\z'_{\sigma(1)}, \dots, \z'_{\sigma(k)})}$
    weakly converges as $n' \to \infty$,
    to a distribution $\PP_{\z',k,\infty}$ on $\Rb^{k p}$.
    There exists a distribution $\PP_{\z',\infty}$ on $\Rb^{k p}$, with a density $f_{\z', \infty}$ with respect to the $p$-dimensional Lebesgue measure such that $\PP_{\z',k,\infty} = \PP_{\z',\infty}^{\otimes k}.$
    \item The matrix $V_1 := \int \psibm(\z'_1, \dots \z'_k) \psibm(\z'_1, \dots \z'_k)^T f_{\z', \infty}(\z'_1) \cdots f_{\z', \infty}(\z'_k) d\z'_1 \cdots d\z'_k$ is non-singular.
    \item \label{condUstat:assumpt:Lambda_C2}
    $\Lambda(\cdot)$ is two times continuously differentiable. Let $\Tc$ be the range of $\theta$, from $\Zc^k$ towards $\Rb$.
    On an open neighborhood of $\Tc$, the second derivative of $\Lambda(\cdot)$ is bounded by a constant $C_{\Lambda''}$.
    \item \label{condUstat:assumpt:finite_integrals} Several integrals exist and are finite, including
    \begin{align*}
        \tilde V_1 &:= \int \theta \big( \z'_1, \dots, \z'_k \big) 
        \Lambda' \Big( \theta \big( \z'_1, \dots, \z'_k \big) \Big)
        \psibm \big( \z'_1, \dots, \z'_k \big)
        f_{\z', \infty}(\z_1) \cdots f_{\z', \infty}(\z_k) \,
        d\z_1 \cdots d\z_k \, \text{ and } \\
        V_2 &:= \int \dfrac{||K||_2^2}{f_\Z(\z'_1)}
        g \big(\x_1, \x_2, \dots, \x_{k} \big)
        g \big(\x_1, \y_2, \dots, \y_{k} \big)
        \Lambda'^2 \big( \theta (\z'_1, \dots, \z'_k) \big)
        \psibm \big(\z'_1, \dots, \z'_k \big)
        \psibm \big(\z'_1, \dots, \z'_k \big)^T \\
        & \hspace{1.5cm} \times f_{\X | \Z = \z'_1}(\x_1)
        \, d\mu(\x_1) d\mu(\z'_1) \prod_{i=2}^k
        f_{\X | \Z = \z'_i}(\y_i) f_{\X | \Z = \z'_i}(\x_i)
        f_{\z', \infty}(\z'_i) \, d\mu(\x_i) d\mu(\y_i) d\z'_{i}.
    \end{align*}
\end{enumerate}
\label{condUstat:assumpt:asymptNorm_joint}
\end{assumpt}

Define $\tilde r_{n,n'}
:= (n \times n' \times h_{n,n'}^{p})^{1/2}$,
$\u := \tilde r_{n,n'} (\beta - \beta^*)$ and
$\hat \u_{n,n'} := \tilde r_{n,n'} (\hat \beta_{n, n'} - \beta^*)$,
so that $\hat \beta_{n, n'} = \beta^* + \hat\u_{n,n'} / \tilde r_{n,n'}$.
We define for every $\u \in \Rb^{p'}$,
\begin{align}
    \FF_{n,n'}(\u)
    & := \frac{- 2 \tilde r_{n,n'}}{|\Iknpr|}
    \sum_{\sigma \in \Iknpr} \xi_{\sigma,n} \psibm \big(\z'_{\sigma(1)}, \dots, \z'_{\sigma(k)} \big)^T \u \nonumber \\
    & + \frac{1}{|\Iknpr|} \sum_{\sigma \in \Iknpr}
    \big\{ \psibm \big(\z'_{\sigma(1)}, \dots, \z'_{\sigma(k)} \big)^T \u \big\}^2
    + \lambda_{n, n'} \tilde r_{n,n'}^2
    \Big( \big| \beta^* + \frac{\u}{\tilde r_{n,n'} } \big|_1 - \big| \beta^* \big|_1 \Big),
    \label{condUstat:eq:def_F_n_nprime}
\end{align}
and we obtain
$\hat\u_{n,n'} = \arg \min_{\u \in \Rb^{p'}} \FF_{n,n'}(\u)$ applying Lemma \ref{condUstat:def:process_GG}.

\begin{lemma}
    Under the same assumptions as in Theorem \ref{condUstat:thm:weak_conv_doubleAsympt},
    \begin{equation*}
        T_1 := \frac{\tilde r_{n,n'}}{|\Iknpr|}
        \sum_{\sigma \in \Iknpr} \xi_{\sigma,n}
        \psibm \big(\z'_{\sigma(1)}, \dots, \z'_{\sigma(k)} \big)
        \indistrto \Nc(0, V_2).
    \end{equation*}
    \label{lemma:limit_T_1}
\end{lemma}
This lemma is proved in~\ref{condUstat:proof:lemma:limit_T_1}.
It will help to control the first term of Equation (\ref{condUstat:eq:def_F_n_nprime}), which is simply $-2 T_1^T \u$.

\mds

Concerning the second term of Equation (\ref{condUstat:eq:def_F_n_nprime}),
using Assumption \ref{condUstat:assumpt:asymptNorm_joint}(iii), we have for every $\u \in \Rb^{p'}$
\begin{align}
    \frac{1}{|\Iknpr|} \sum_{\sigma \in \Iknpr}
    &\big\{ \psibm \big(\z'_{\sigma(1)}, \dots, \z'_{\sigma(k)} \big)^T \u \big\}^2
    \nonumber \\
    &\to \int \left( \psibm (\z'_1, \dots, \z'_k)^T \u \right)^2
    f_{\z', \infty}(\z'_1) \cdots f_{\z', \infty}(\z'_k) \,
    d\z'_1 \cdots d\z'_k.
    \label{condUstat:eq:limit_2rdterm_F}
\end{align}
This has to be read as a convergence of a sequence of real numbers indexed by $\u$, because the design points $\z'_i$ are deterministic.
We also have, for any $\u \in \Rb^{p'}$ and when $n$ is large enough,
\begin{align*}
    \big| \beta^* + \frac{\u}{\tilde r_{n,n'} } \big|_1 - \big|  \beta^* \big|_1
    = \sum_{i=1}^{p'} \Big( \frac{|u_i|}{\tilde r_{n,n'}} \1_{\{\beta^*_i = 0\}}
    + \frac{u_i}{\tilde r_{n,n'}}  \sgn(\beta^*_i) \1_{\{\beta^*_i \neq 0\}} \Big).
\end{align*}

Therefore, by Assumption \ref{condUstat:assumpt:asymptNorm_joint}(ii)(b), for every $\u \in \Rb^{p'}$,
\begin{align}
    \lambda_{n, n'} \tilde r_{n,n'}^2
    \Big( \big| \beta^* + \frac{\u}{\tilde r_{n,n'} } \big|_1 - \big| \beta^* \big|_1 \Big) \to 0,
    \label{condUstat:eq:limit_3rdterm_F}
\end{align}
when $(n,n')$ tends to the infinity.
Combining Lemma \ref{lemma:limit_T_1} and Equations (\ref{condUstat:eq:def_F_n_nprime}-\ref{condUstat:eq:limit_3rdterm_F}), and defining the function $\FF_{\infty, \infty}$ by
\begin{align*}
    \FF_{\infty, \infty}(\u)
    := 2 \tilde \W^T \u
    + \int \left( \psibm (\z'_1, \dots, \z'_k)^T \u \right)^2
    f_{\z', \infty}(\z'_1) \cdots f_{\z', \infty}(\z'_k) \,
    d\z'_1 \cdots d\z'_k,
\end{align*}
where $\u \in \Rb^{r}$ and $\tilde W \sim \Nc(0, V_2)$, we obtain that every finite-dimensional margin of $\FF_{n,n'}$ weakly converges to the corresponding margin of $\FF_{\infty, \infty}$.
Now, applying the convexity lemma, we get
\begin{align*}
    \hat \u_{n,n'} \indistrto \u_{\infty, \infty}, \text{ where }
    \u_{\infty, \infty} := \arg \min_{\u \in \Rb^{r}} \FF_{\infty, \infty}(\u).
\end{align*}
Since $\FF_{\infty, \infty}(\u)$ is a continuously differentiable convex function, apply the first-order condition $\nabla \FF_{\infty, \infty}(\u) = 0$, which yields
$$2 \tilde \W  + 2 \int \psibm (\z'_1, \dots, \z'_k) \psibm (\z'_1, \dots, \z'_k)^T \u_{\infty, \infty} f_{\z', \infty}(\z'_1) \cdots f_{\z', \infty}(\z'_k) \,
d\z'_1 \cdots d\z'_k = 0.$$
As a consequence $\u_{\infty, \infty} = - V_1^{-1} \tilde \W \sim \Nc(0, \tilde V_{as}),$ using Assumption \ref{condUstat:assumpt:asymptNorm_joint}(iv).
We finally obtain $\tilde r_{n,n'}  \big( \hat \beta_{n, n'} - \beta^* \big) \indistrto \Nc \big(0, \tilde V_{as} \big)$,
as claimed.$\;\; \Box$

\subsection{Proof of Lemma \ref{lemma:limit_T_1}}
\label{condUstat:proof:lemma:limit_T_1}

Using a Taylor expansion yields
\begin{align*}
    T_1  &:= \frac{\tilde r_{n,n'}}{|\Iknpr|}
    \sum_{\sigma \in \Iknpr} \xi_{\sigma,n}
    \psibm \big(\z'_{\sigma(1)}, \dots, \z'_{\sigma(k)} \big) \\
    &= \frac{\tilde r_{n,n'}}{|\Iknpr|}
    \sum_{\sigma \in \Iknpr} \bigg( \Lambda \Big( \hat \theta
    \big(\z'_{\sigma(1)}, \dots, \z'_{\sigma(k)} \big) \Big)
    - \Lambda \Big( \theta
    \big(\z'_{\sigma(1)}, \dots, \z'_{\sigma(k)} \big) \Big) \bigg)
    \psibm \big(\z'_{\sigma(1)}, \dots, \z'_{\sigma(k)} \big) \\
    &= T_2 + T_3,
\end{align*}
where the main term is
\begin{align*}
    T_2 := \frac{\tilde r_{n,n'}}{|\Iknpr|}
    \sum_{\sigma \in \Iknpr} \Lambda' \Big( \theta
    \big(\z'_{\sigma(1)}, \dots, \z'_{\sigma(k)} \big) \Big)
    \Big( \hat \theta \big(\z'_{\sigma(1)}, \dots, \z'_{\sigma(k)} \big)
    - \theta \big(\z'_{\sigma(1)}, \dots, \z'_{\sigma(k)} \big) \Big)
    \psibm \big(\z'_{\sigma(1)}, \dots, \z'_{\sigma(k)} \big),
\end{align*}
and the remainder is
\begin{align*}
    T_3 := \frac{\tilde r_{n,n'}}{|\Iknpr|}
    \sum_{\sigma \in \Iknpr} \alpha_{3, \sigma} \cdot
    \Big( \hat \theta \big(\z'_{\sigma(1)}, \dots, \z'_{\sigma(k)} \big)
    - \theta \big(\z'_{\sigma(1)}, \dots, \z'_{\sigma(k)} \big) \Big)^2
    \psibm \big(\z'_{\sigma(1)}, \dots, \z'_{\sigma(k)} \big),
\end{align*}
with $\forall \sigma \in \Iknpr$, $|\alpha_{3,\sigma}| \leq C_{\Lambda''}/2$, by Assumption \ref{condUstat:assumpt:asymptNorm_joint}\ref{condUstat:assumpt:Lambda_C2}.

\mds

Let us define
$\bar \psibm_\sigma := \Lambda' \Big( \theta
\big(\z'_{\sigma(1)}, \dots, \z'_{\sigma(k)} \big) \Big)
\psibm \big(\z'_{\sigma(1)}, \dots, \z'_{\sigma(k)} \big),$
for every $\sigma \in \Iknpr$.
Using the definition (\ref{condUstat:def:estimator_hat_theta}), we rewrite $T_2 := T_4 + T_5$ where
\begin{align*}
    T_4 &:= \frac{\tilde r_{n,n'}}{|\Iknpr|\cdot|\Ikn|}
    \sum_{\sigma \in \Iknpr} \sum_{\varsigma \in \Ikn}
    \dfrac{\prod_{i=1}^k K_h \big(\Z_{\varsigma(i)}-\z'_{\sigma(i)} \big) }
    { \prod_{i=1}^k f_\Z(\z'_{\sigma(i)}) }
    \Big( g \big(\X_{\varsigma(1)}, \dots, \X_{\varsigma(k)} \big) - \theta \big(\z'_{\sigma(1)}, \dots, \z'_{\sigma(k)} \big) \Big)
    \bar \psibm_\sigma, \\
    T_5 &:= \frac{\tilde r_{n,n'}}{|\Iknpr|\cdot|\Ikn|}
    \sum_{\sigma \in \Iknpr} \sum_{\varsigma \in \Ikn}
    \prod_{i=1}^k K_h \big(\Z_{\varsigma(i)}-\z'_{\sigma(i)} \big)
    \Big( g \big(\X_{\varsigma(1)}, \dots, \X_{\varsigma(k)} \big) - \theta \big(\z'_{\sigma(1)}, \dots, \z'_{\sigma(k)} \big) \Big) \\
    &  \hspace{3cm} \times \bigg( \dfrac{1} {N_k(\z'_{\sigma(1)}, \dots, \z'_{\sigma(k)})}
    - \dfrac{1} {\prod_{i=1}^k f_\Z(\z'_{\sigma(i)}) }\bigg)
    \bar \psibm_\sigma.
\end{align*}
To lighten the notations, we will define
$K_{\sigma, \varsigma}
:= \prod_{i=1}^k K_h \big(\Z_{\varsigma(i)}-\z'_{\sigma(i)} \big),$ 
$g_\varsigma := g \big(\X_{\varsigma(1)}, \dots, \X_{\varsigma(k)} \big),$
$\theta_\sigma := \theta \big(\z'_{\sigma(1)}, \dots, \z'_{\sigma(k)} \big),$
$f_{\Z', \sigma} := \prod_{i=1}^k f_\Z(\z'_{\sigma(i)}),$
and $N_\sigma := N_k(\z'_{\sigma(1)}, \dots, \z'_{\sigma(k)})$,
for every $\sigma \in \Iknpr$ and $\varsigma \in \Ikn$, so that
\begin{align}
    T_4 &:= \frac{\tilde r_{n,n'}}{|\Iknpr|\cdot|\Ikn|}
    \sum_{\sigma \in \Iknpr} \sum_{\varsigma \in \Ikn}
    \dfrac{ K_{\sigma, \varsigma} } { f_{\Z', \sigma} }
    \big( g_\varsigma - \theta_\sigma \big)
    \bar \psibm_\sigma, \\
    T_5 &:= \frac{\tilde r_{n,n'}}{|\Iknpr|\cdot|\Ikn|}
    \sum_{\sigma \in \Iknpr} \sum_{\varsigma \in \Ikn}
    K_{\sigma, \varsigma}
    \big( g_\varsigma - \theta_\sigma \big)
    \Big( \dfrac{1} {N_\sigma}
    - \dfrac{1} {f_{\Z', \sigma}} \Big)
    \bar \psibm_\sigma. \label{condUstat:eq:def:term_T_5}
\end{align}

Using $\alpha$-order limited expansions, we get
\begin{align}
    \EE[T_4]
    &= \frac{\tilde r_{n,n'}}{|\Iknpr|} \sum_{\sigma \in \Iknpr}
    \int \frac{\prod_{i=1}^k K_h \big(\z_i-\z'_{\sigma(i)} \big)} {f_{\Z', \sigma}}
    \Big( g \big(\x_{1:k} \big) - \theta_\sigma \Big) 
    \prod_{i=1}^k f_{\X, \Z} (\x_i, \z_i)
    d\mu^{\otimes k} (\x_{1:k}) d\z_{1:k} \label{condUstat:eq:expectancy_T_4} \\
    &= \frac{\tilde r_{n,n'}}{|\Iknpr|} \sum_{\sigma \in \Iknpr}
    \int \frac{\prod_{i=1}^k K \big(\t_i\big)} {f_{\Z', \sigma}}
    \Big( g \big(\x_{1:k} \big) - \theta_\sigma \Big) 
    \prod_{i=1}^k f_{\X, \Z} (\x_i, \z'_{\sigma(i)} + h \t_i)
    d\mu^{\otimes k} (\x_{1:k}) d\t_{1:k} \nonumber \\
    &= \frac{\tilde r_{n,n'} h^{k \alpha}}{|\Iknpr|} \sum_{\sigma \in \Iknpr}
    \int \frac{\prod_{i=1}^k K \big(\t_i\big)} {f_{\Z', \sigma}}
    \Big( g \big(\x_{1:k} \big) - \theta_\sigma \Big) 
    \prod_{i=1}^k d^{(\alpha)}_\Z f_{\X, \Z} (\x_i, \z^*_{\sigma(i)})
    d\mu^{\otimes k} (\x_{1:k}) d\t_{1:k} \nonumber \\
    &= O \Big(  \tilde r_{n,n'} h^{k \alpha} \Big)
    = O \Big( (n \times n' \times h_{n,n'}^{p + 2 k \alpha})^{1/2} \Big) = o(1), \nonumber
\end{align}
where above, $\z_i^*$ denote some vectors in $\Rb^p$ such that $|\z'_i - \z^*_i|_\infty \leq 1$, depending on $\z'_i$ and $\x_i$.

\mds

We can therefore use the centered version of $T_4$, defined as
\begin{align*}
    T_4 - \EE[T_4] &= \frac{\tilde r_{n,n'}}{|\Iknpr|\cdot|\Ikn|}
    \sum_{\sigma \in \Iknpr} \sum_{\varsigma \in \Ikn}
    g_{\sigma, \varsigma}, \\
    g_{\sigma, \varsigma}
    &:= \dfrac{\bar \psibm_\sigma} { f_{\Z', \sigma} }
    \Big( K_{\sigma, \varsigma} \big( g_\varsigma - \theta_\sigma \big)
    - \EE \big[ K_{\sigma, \varsigma} \big( g_\varsigma - \theta_\sigma \big) \big] \Big).
\end{align*}

\mds

{\bf Computation of the limit of the variance matrix $Var[T_4]$.}

\mds 

We have $Var[T_4] = \EE[T_4 T_4^T] + o(1)$.
\begin{align*}
    Var[T_4] = \frac{\tilde r_{n,n'}^2}{|\Iknpr|^2\cdot|\Ikn|^2}
    \sum_{\sigma, \bar \sigma \in \Iknpr} \sum_{\varsigma, \bar \varsigma \in \Ikn}
    \EE[ g_{\sigma, \varsigma}
    g_{\bar \sigma, \bar \varsigma}^T ]
    + o(1).
\end{align*}
By independence, $\EE[ g_{\sigma, \varsigma}
g_{\bar \sigma, \bar \varsigma}^T ] = 0$ as soon as $\varsigma \cap \bar \varsigma = \emptyset$, where we identify a permutation $\varsigma$ and its image  $\varsigma(\{1, \dots, k\})$. Therefore, we get
\begin{align*}
    Var[T_4]
    &\simeq \frac{n n' h_{n,n'}^{p}} {|\Iknpr|^2 \cdot |\Ikn|^2}
    \sum_{\sigma, \bar \sigma \in \Iknpr}
    \sum_{\substack{ \varsigma, \bar \varsigma \in \Ikn \\
    \varsigma \cap \bar \varsigma \neq \emptyset}}
    \EE \big[ g_{\sigma, \varsigma} \, 
    g_{\bar \sigma, \bar \varsigma}^T \big] \\
    &= \frac{n n' h_{n,n'}^{p}} {|\Iknpr|^2 \cdot |\Ikn|^2}
    \sum_{\sigma, \bar \sigma \in \Iknpr}
    \sum_{\substack{ \varsigma, \bar \varsigma \in \Ikn \\
    \varsigma \cap \bar \varsigma \neq \emptyset}}
    g_{\sigma, \varsigma, \bar \sigma, \bar \varsigma}
    - \tilde g_\sigma \tilde g_{\bar \sigma}^T,
\end{align*}
where $\tilde g_\sigma := \bar \psibm_\sigma
\EE \big[ K_{\sigma, \varsigma} \big( g_\varsigma - \theta_\sigma \big) \big] / f_{\Z', \sigma}$ and
\begin{align*}
    g_{\sigma, \varsigma, \bar \sigma, \bar \varsigma}
    := \dfrac{\bar \psibm_\sigma \bar \psibm_{\bar \sigma}^T}
    { f_{\Z', \sigma} f_{\Z', \bar \sigma} }
    \EE \bigg[ K_{\sigma, \varsigma}
    K_{\bar \sigma, \bar \varsigma}
    \big( g_\varsigma - \theta_\sigma \big)
    \big( g_{\bar \varsigma} - \theta_{\bar \sigma} \big)
    \bigg].
\end{align*}
Assume now that $\varsigma \cap \bar \varsigma$ is of cardinality $1$, i.e. there exists only one couple $(j, \bar j) \in \{1, \dots, k\}^2$ such that $\varsigma(j) = \bar \varsigma(\bar j)$.
Then,
\begin{align*}
    g_{\sigma, \varsigma, \bar \sigma, \bar \varsigma}
    &= \dfrac{\bar \psibm_\sigma \bar \psibm_{\bar \sigma}^T}
    { f_{\Z', \sigma} f_{\Z', \bar \sigma} }
    \int \big( g(\X_{1:k}) - \theta_\sigma \big)
    \big( g(\x_{k+1}, \dots, \x_{k + \bar j - 1}, \x_j,
    \x_{k + \bar j + 1}, \dots, \x_{2k}) - \theta_{\bar \sigma} \big) \\
    &\hspace{2cm} \cdot \prod_{i=1}^k K_h \big(\z_{i}-\z'_{\sigma(i)} \big)
    f_{\X, \Z} (\x_i, \z_i) d\mu(\x_i) d\z_i
    \cdot K_h \big(\z_{j}-\z'_{\bar \sigma(\bar j)} \big)\\
    &\hspace{2cm} \cdot \prod_{\bar i=1, \, \bar i \neq \bar j}^{k}
    K_h \big(\z_{k+i}-\z'_{\bar \sigma(\bar i)} \big)
    f_{\X, \Z} (\x_{k+i}, \z_{k+i}) d\mu(\x_{k+i}) d\z_{k+i} \\
    &= \dfrac{\bar \psibm_\sigma \bar \psibm_{\bar \sigma}^T}
    { f_\Z(\z_j) }
    \int \big( g(\X_{1:k}) - \theta_\sigma \big)
    \big( g(\x_{k+1}, \dots, \x_{k + \bar j - 1}, \x_j,
    \x_{k + \bar j + 1}, \dots, \x_{2k}) - \theta_{\bar \sigma} \big) \\
    &\hspace{2cm} \cdot \prod_{i=1}^k K(\t_{i})
    \frac{f_{\X, \Z} (\x_i, \z'_{\sigma(i)} + h\t_i)}{f_\Z(\z'_{\sigma(i)})} d\mu(\x_i) d\t_i
    \cdot h^{-p} K \bigg(\t_i +
    \frac{\z'_{\sigma(j)} - \z'_{\bar \sigma(\bar j)}} {h} \bigg)\\
    &\hspace{2cm} \cdot \prod_{\bar i=1, \, \bar i \neq \bar j}^{k}
    K(\t_{k+i})
    \frac{f_{\X, \Z} (\x_{k+i}, \z'_{\bar \sigma(\bar i)} + h \t_{k+i})} {f_\Z(\z_{k+i})}
    d\mu(\x_{k+i}) d\t_{k+i} \\
    &\simeq \dfrac{\bar \psibm_\sigma \bar \psibm_{\bar \sigma}^T}
    { f_\Z(\z_j) }
    \int \big( g(\X_{1:k}) - \theta_\sigma \big)
    \big( g(\x_{k+1}, \dots, \x_{k + \bar j - 1}, \x_j,
    \x_{k + \bar j + 1}, \dots, \x_{2k}) - \theta_{\bar \sigma} \big) \\
    &\hspace{2cm} \cdot \prod_{i=1}^k K(\t_{i})
    \frac{f_{\X, \Z} (\x_i, \z'_{\sigma(i)})}{f_\Z(\z_i)} d\mu(\x_i) d\t_i
    \cdot h^{-p} K \bigg(\t_i +
    \frac{\z'_{\sigma(j)} - \z'_{\bar \sigma(\bar j)}} {h} \bigg) \\
    &\hspace{2cm} \cdot \prod_{\bar i=1, \, \bar i \neq \bar j}^{k}
    K(\t_{k+i})
    \frac{f_{\X, \Z} (\x_{k+i}, \z'_{\bar \sigma(\bar i)})} {f_\Z(\z'_{\bar \sigma(\bar i)})}
    d\mu(\x_{k+i}) d\t_{k+i}.
\end{align*}
By assumption, this is zero unless $\sigma(j) = \bar \sigma(\bar j)$. In this case, it can be simplified, giving
\begin{align*}
    g_{\sigma, \varsigma, \bar \sigma, \bar \varsigma}
    &\simeq \dfrac{\bar \psibm_\sigma \bar \psibm_{\bar \sigma}^T}
    { f_\Z(\z_j) h^{p} } \int K^2
    \int \big( g(\x_{1:k}) - \theta_\sigma \big)
    \big( g( \x_{k:2k, \bar j \to j} ) - \theta_{\bar \sigma} \big) \\
    &\cdot \prod_{i=1}^k
    f_{\X | \Z = \z'_{\sigma(i)}} (\x_{k}) d\mu(\x_i)
    \prod_{\bar i=1, \, \bar i \neq \bar j}^{k}
    f_{\X | \Z = \z'_{\bar \sigma(\bar i)}} (\x_{k+i})
    d\mu(\x_{k+i})
    \, =: h^{-p} g_{\sigma, \bar \sigma, j, \bar j},
\end{align*}
where $\x_{k:2k, \bar j \to j}
:= (\x_{k+1}, \dots, \x_{k + \bar j - 1}, \x_j,
\x_{k + \bar j + 1}, \dots, \x_{2k})$.

\mds

Note that, if $\varsigma \cap \bar \varsigma$ is of cardinality strictly greater than $1$, some supplementary powers of $h^{-p}$ arise thanks to the repeated kernels in $\varsigma$ and $\bar \varsigma$. As a consequence, they are of lower order and therefore negligible.
Using $\alpha$-order expansions as in Equation (\ref{condUstat:eq:expectancy_T_4}), we get $\sup_\sigma |\tilde g_\sigma| = O(h^{k \alpha})$. Thus,

\begin{align*}
    Var[T_4]
    &\simeq O \big( n n' h_{n,n'}^{p + 2 k \alpha} \big)
    + \frac{n n' h_{n,n'}^{p}} {|\Iknpr|^2 \cdot |\Ikn|^2}
    \sum_{\varsigma \in \Ikn} \sum_{j, \bar j = 1}^k
    \sum_{\substack{\bar \varsigma \in \Ikn \\
    \varsigma(j) = \bar \varsigma(\bar j), \,
    | \varsigma \cap \bar \varsigma | = 1} }
    \sum_{\sigma, \bar \sigma \in \Iknpr, \,
    \sigma(j) = \bar \sigma(\bar j)}
    h^{-p} g_{\sigma, \bar \sigma, j, \bar j}
    \displaybreak[0] \\
    &\simeq \frac{n'} {|\Iknpr|^2}
    \sum_{j, \bar j = 1}^k
    \sum_{\sigma, \bar \sigma \in \Iknpr, \,
    \sigma(j) = \bar \sigma(\bar j)}
    g_{\sigma, \bar \sigma, j, \bar j} \\
    &\to \sum_{j, \bar j = 1}^k
    g_{j, \bar j,\infty} = V_2,
\end{align*}
where
\begin{align*}
    g_{j, \bar j,\infty}
    &:= \int \Lambda' \Big( \theta \big(\z'_{1:k} \big) \Big)
    \Lambda' \Big( \theta \big(\z'_{k:2k, \bar j \to j} \big) \Big) \psibm \big(\z'_{1:k} \big)
    \psibm^T \big(\z'_{k:2k, \bar j \to j} \big)
    \dfrac{\int K^2} { f_\Z(\z'_j) }
    \int \big( g(\x_{1:k}) - \theta (\z'_{1:k}) \big) \\
    & \cdot \big( g(\x_{k:2k, \bar j \to j})
    - \theta (\z'_{k:2k, \bar j \to j}) \big) 
    \prod_{i=1, i \neq k +\bar j}^{2k}
    f_{\X | \Z = \z'_{i}} (\x_{i}) f_{\Z', \infty}(\z'_i) d\mu(\x_i) d\z'_i.
\end{align*}
In Section \ref{section:proof:asympt_norm_T_4}, we will prove that $T_4$ is asymptotically Gaussian ; therefore, its asymptotic variance will be given by $V_2$.

\mds

Now, decompose the term $T_5$, defined in Equation~(\ref{condUstat:eq:def:term_T_5}), using a Taylor expansion of the function $x \mapsto 1/(1+x)$ at $0$.
\begin{align*}
    \dfrac{1} {N_\sigma} - \dfrac{1} {f_{\Z', \sigma}}
    = \dfrac{1} {f_{\Z', \sigma}} \Bigg(
    \dfrac{1}{1 + \frac{N_\sigma - f_{\Z', \sigma}}{f_{\Z', \sigma}}}
    - 1 \Bigg)
    = - \dfrac{N_\sigma - f_{\Z', \sigma}}{f_{\Z', \sigma}^2}
    + T_{7, \sigma},
\end{align*}
where 
\begin{align*}
    T_{7, \sigma} := \frac{1}{f_{\Z', \sigma}} (1 + \alpha_{7,\sigma})^{-3}
    \left( \dfrac{N_\sigma - f_{\Z', \sigma}}{f_{\Z', \sigma}}
    \right)^2,
    \text{ with } |\alpha_{7,\sigma}| \leq
    \left| \dfrac{N_\sigma - f_{\Z', \sigma}}{f_{\Z', \sigma}} \right|.
\end{align*}

We have therefore the decomposition $T_5 = - T_6 + T_7$, where
\begin{align}
    T_6 &:= \frac{\tilde r_{n,n'}}{|\Iknpr|\cdot|\Ikn|}
    \sum_{\sigma \in \Iknpr} \sum_{\varsigma \in \Ikn}
    K_{\sigma, \varsigma}
    \big( g_\varsigma - \theta_\sigma \big)
    \dfrac{N_\sigma - f_{\Z', \sigma}}{f_{\Z', \sigma}^2}
    \bar \psibm_\sigma , \label{condUstat:eq:def:term_T_6} \\
    T_7 &:= \frac{\tilde r_{n,n'}}{|\Iknpr|\cdot|\Ikn|}
    \sum_{\sigma \in \Iknpr} \sum_{\varsigma \in \Ikn}
    K_{\sigma, \varsigma}
    \big( g_\varsigma - \theta_\sigma \big) T_{7,\sigma}
    \bar \psibm_\sigma. \label{condUstat:eq:def:term_T_7} 
\end{align}

Summing up all the previous equation, we get
\begin{align*}
    T_1 = (T_4 - \EE[T_4]) - T_6 + T_7 + T_3 + o(1).
\end{align*}
Afterwards, we will prove that all the remainders terms $T_6$, $T_7$ and $T_3$ are negligible, i.e. they tend to zero in probability.
These results are respectively proved in Subsections \ref{section:proof:limit_T_6}, \ref{section:proof:limit_T_7} and \ref{section:proof:limit_T_3}.
Combining all these elements with the asymptotic normality of $T_4$ (proved in Subsection~\ref{section:proof:asympt_norm_T_4}), we get $T_1 \indistrto \Nc(0, V_2)$, as claimed.$\;\; \Box$

\subsection{Proof of the asymptotic normality of \texorpdfstring{$T_4$}{T4}}
\label{section:proof:asympt_norm_T_4}

Using the H\'ajek projection of $T_4$, we define
\begingroup \allowdisplaybreaks
\begin{align*}
    &T_4 - \EE[T_4]
    = T_{4,1} + T_{4,2}, \text{ where } \\
    &T_{4,1} := \frac{\tilde r_{n,n'}}{|\Iknpr|\cdot|\Ikn|}
    \sum_{\sigma \in \Iknpr} \sum_{\varsigma \in \Ikn}
    \sum_{i=1}^k
    \EE[ g_{\sigma, \varsigma} | \varsigma(i)] , \\
    &T_{4,2} := \frac{\tilde r_{n,n'}}{|\Iknpr|\cdot|\Ikn|}
    \sum_{\sigma \in \Iknpr} \sum_{\varsigma \in \Ikn}
    \bigg( g_{\sigma, \varsigma} -
    \sum_{i=1, \dots, k} \EE[ g_{\sigma, \varsigma} | \varsigma(i)] 
    \bigg),
\end{align*}
\endgroup
denoting by $|i$ the conditioning with respect to $(\X_i, \Z_i)$, for $i \in \{1, \dots, n \}$.
We will show that $T_{4,1}$ is asymptotically normal, and that $T_{4,2}=o(1)$.

\mds

Using the fact that the $(\X_i, \Z_i)_i$ are i.i.d., and denoting by~$Id$ the injective function $i \mapsto i$, we have
\begin{align*}
    T_{4,1}
    &= \frac{k \tilde r_{n,n'}}{n |\Iknpr|}
    \sum_{\sigma \in \Iknpr} \sum_{i=1}^n
    \EE \bigg[ \dfrac{\bar \psibm_\sigma} { f_{\Z', \sigma} }
    K_{\sigma, Id} \big( g_{Id} - \theta_\sigma \big)
    - \bar g_\sigma \bigg| i \bigg] \\
    &\simeq \frac{k \tilde r_{n,n'}}{n |\Iknpr|}
    \sum_{\sigma \in \Iknpr} \sum_{i=1}^n
    \EE \bigg[ \dfrac{\bar \psibm_\sigma} { f_{\Z', \sigma} }
    K_{\sigma, Id} \big( g_{Id} - \theta_\sigma \big)
    \bigg| i \bigg] =: \sum_{i=1}^n \alpha_{4,i,n},
\end{align*}
because $\sup_\sigma |\bar g_\sigma| = O(h^{k \alpha})$, as proved in the previous section, hence negligible.
The $\alpha_{4,i,n},$ for $1 \leq i \leq n$, form a triangular array of i.i.d. variables.
To prove the asymptotic normality of $T_{4,1}$, it remains to check Lyapunov's condition, i.e. we will show that
$\sum_{i=1}^n \EE \big[ |\alpha_{4,i,n}|_\infty^3 \big] \to 0$.
We have
\begingroup \allowdisplaybreaks
\begin{align*}
    \sum_{i=1}^n & \, \EE \big[ |\alpha_{4,i,n}|_\infty^3 \big]
    = n \, \EE \big[ |\alpha_{4,1,n}|_\infty^3 \big] \\
    &= \frac{k^3 n \tilde r_{n,n'}^3}{n^3 |\Iknpr|^3}
    \sum_{\sigma, \nu, \vartheta \in \Iknpr}
    \dfrac{\bar \psibm_\sigma \otimes \bar \psibm_\nu \otimes \bar \psibm_\vartheta} { f_{\Z', \sigma} f_{\Z', \nu} f_{\Z', \vartheta} }
    \EE \Bigg[ 
    \EE \Big[ K_{\sigma, Id} 
    \big( g_{Id} - \theta_\sigma \big) \Big| 1 \Big]
    \EE \Big[ K_{\nu, Id}
    \big(  g_{Id} - \theta_\nu) \big) \Big| 1 \Big]
    \EE \Big[ K_{\vartheta, Id}
    \big( g_{Id} - \theta_\vartheta) \big) \Big| 1 \Big] \Bigg] \\
    &= \frac{k^3 \tilde r_{n,n'}^3}{n^2 |\Iknpr|^3}
    \sum_{\sigma, \nu, \vartheta \in \Iknpr}
    \dfrac{\bar \psibm_\sigma \otimes \bar \psibm_\nu \otimes \bar \psibm_\vartheta} { f_{\Z}(\z'_{\nu(1)}) f_{\Z}(\z'_{\vartheta(1)}) }
    \int K_h \big(\z_1-\z'_{\sigma(1)} \big)
    K_h \big(\z_1-\z'_{\nu(1)} \big)
    K_h \big(\z_1-\z'_{\vartheta(1)} \big) \\
    &\cdot \prod_{i=2}^{k} K_h \big(\z_{i}-\z'_{\sigma(i)} \big)
    K_h \big(\z_{k+i}-\z'_{\nu(i)} \big)
    K_h \big(\z_{2k+i}-\z'_{\vartheta(i)} \big) \\
    &\cdot \big( g(\x_{1:k}) - \theta_\sigma) \big)
    \big( g(\x_1, \x_{(k+2):(2k)}) - \theta_\nu)
    \big( g(\x_1, \x_{(2k+2):(3k)}) - \theta_\vartheta) \\
    & \cdot \prod_{i=1}^{k} \frac{f_{\X, \Z}(\x_i, \z_i)}
    {f_\Z(\z'_{\sigma(i)})} d\mu(\x_i) d\z_i
    \prod_{i=2}^{k} \frac{f_{\X, \Z}(\x_{k+i}, \z_{k+i})}
    {f_\Z(\z'_{\nu(i)})} d\mu(\x_{k+i}) d\z_{k+i}
    \prod_{i=2}^{k} \frac{f_{\X, \Z}(\x_{2k+i}, \z_{2k+i})}
    {f_\Z(\z'_{\vartheta(i)})} d\mu(\x_{2k+i}) d\z_{2k+i} \\
    &\simeq \frac{k^3 \tilde r_{n,n'}^3}{n^2 |\Iknpr|^3}
    \sum_{\sigma, \nu, \vartheta \in \Iknpr}
    \dfrac{\bar \psibm_\sigma \otimes \bar \psibm_\nu \otimes \bar \psibm_\vartheta} { f_{\Z}(\z'_{\nu(1)}) f_{\Z}(\z'_{\vartheta(1)}) }
    \int h^{-2 p} K( \t_1 )
    K \left( \t_1 + \frac{\z'_{\sigma(1)} - \z'_{\nu(1)}}{h} \right)
    K \left( \t_1 + \frac{\z'_{\sigma(1)} - \z'_{\vartheta(1)}}{h} \right) \\
    &\cdot \prod_{i=2}^{k} K_h \big( \t_{i} \big) K_h \big(\t_{k+i} \big)
    K_h \big(\t_{2k+i} \big)
    \big( g(\x_{1:k}) - \theta_\sigma) \big)
    \big( g(\x_1, \x_{(k+2):(2k)}) - \theta_\nu)
    \big( g(\x_1, \x_{(2k+2):(3k)}) - \theta_\vartheta) \\
    & \cdot \prod_{i=1}^{k} f_{\X | \Z = \z'_{\sigma(i)}}(\x_i) d\mu(\x_i) d\z_i
    \prod_{i=2}^{k} f_{\X | \Z = \z'_{\nu(i)}}(\x_{k+i}) d\mu(\x_{k+i}) d\z_{k+i}
    \prod_{i=2}^{k} f_{\X | \Z = \z'_{\vartheta(i)}}(\x_{2k+i}, \t_{2k+i})
    d\mu(\x_{2k+i}) d\t_{2k+i},
\end{align*}
\endgroup
where in the last equivalent, we use a change of variable from the $\z_i$ to the $\t_i$, and then the continuity of the density $f_{\X,\Z}$ with respect to $\z$, because $h = o(1)$.

\mds

Because of our assumptions, the terms of the sum for which $\sigma(1) \neq 1$ or $\nu(1) \neq 1$ are zero. Therefore, we get 
\begin{align*}
    \sum_{i=1}^n & \, \EE \big[ |\alpha_{4,i,n}|_\infty^3 \big]
    =  \frac{\tilde r_{n,n'}^3 h^{-2 p}}{n^2 |\Iknpr|^3}
    \sum_{\sigma, \nu, \vartheta \in \Iknpr, \sigma(1) = \nu(1) = 1} O(1)
    = O \left( \frac{(n n' h^{p})^{3/2}}{n^2 n'{}^2 h^{2 p}} \right)
    = O \left( \frac{1}{(n n' h^{p})^{1/2}} \right)
    = o(1).
\end{align*}

\mds

We prove now that $T_{4,2} = o(1)$. Note first that, by construction, $\EE[T_{4,2}] = 0$. Computing its variance, we get
\begin{align}
    \EE \big[T_{4,2} T_{4,2}^T \big]
    &= \EE \Bigg[ \frac{\tilde r_{n,n'}^2}{|\Iknpr|^2 \cdot |\Ikn|^2}
    \sum_{\sigma, \bar \sigma \in \Iknpr}
    \sum_{\varsigma, \bar \varsigma \in \Ikn}
    \bigg( g_{\sigma, \varsigma} - \sum_{i=1, \dots, k}
    \EE \big[ g_{\sigma, \varsigma} \big| \varsigma(i) \big] 
    \bigg) \bigg( g_{\bar \sigma, \bar \varsigma} -
    \sum_{\bar i=1, \dots, k} \EE \big[ g_{\bar \sigma, \bar \varsigma}
    \big| \bar \varsigma(\bar i) \big] 
    \bigg)^T \Bigg] \nonumber \\
    &=: \frac{\tilde r_{n,n'}^2}{|\Iknpr|^2 \cdot |\Ikn|^2}
    \sum_{\sigma, \bar \sigma \in \Iknpr}
    \sum_{\varsigma, \bar \varsigma \in \Ikn}
    \EE \Big[ \tilde g_{\sigma, \bar \sigma, \varsigma, \bar \varsigma}
    \Big] \label{condUstat:eq:variance_T_4_2}.
\end{align}
Because of $\EE[g_{\sigma, \varsigma}] = 0$ and by independence, the terms in the latter sum for which $\varsigma \cap \bar \varsigma = \emptyset$ are zero.
Otherwise, there exists $j_1,j_2 \in \{ 1, \dots, k \}$ such that $\varsigma(j_1) = \bar \varsigma(j_2)$.
If $\varsigma \cap \bar \varsigma$ is of cardinal $1$, meaning that there is no other identities between elements of $\varsigma$ and $\bar \varsigma$, then we will show that the corresponding term is zero as well.
We place ourselves in this case, assuming that $|\varsigma \cap \bar \varsigma|=1$, and we get
\begin{align*}
    \EE \Big[ \tilde g_{\sigma, \bar \sigma, \varsigma, \bar \varsigma}
    \Big]
    &= \EE \Bigg[ \bigg( g_{\sigma, \varsigma} -
    \sum_{i=1, \dots, k} \EE \big[ g_{\sigma, \varsigma} \big| \varsigma(i) \big] 
    \bigg) \bigg( g_{\bar \sigma, \bar \varsigma}^T -
    \sum_{\bar i=1, \dots, k} \EE \big[ g_{\bar \sigma, \bar \varsigma}^T \big| \bar \varsigma(\bar i) \big] 
    \bigg) \Bigg] \\
    &= \EE \Bigg[ \bigg( g_{\sigma, \varsigma} -
    \EE \big[ g_{\sigma, \varsigma} \big| \varsigma(j_1) \big] 
    \bigg) \bigg( g_{\bar \sigma, \bar \varsigma}^T -
    \EE \big[ g_{\bar \sigma, \bar \varsigma}^T
    \big| \bar \varsigma(j_2) \big] \bigg) \Bigg] \\
    &= \EE \Bigg[ \EE \bigg[ \Big( g_{\sigma, \varsigma} -
    \EE \big[ g_{\sigma, \varsigma} \big| \varsigma(j_1) \big] 
    \Big) \Big( g_{\bar \sigma, \bar \varsigma}^T -
    \EE \big[ g_{\bar \sigma, \bar \varsigma}^T \big| \varsigma(j_1) \big]
    \Big) \bigg| \varsigma(j_1) \bigg] \Bigg] \\
    &= \EE \Bigg[ \EE \bigg[ g_{\sigma, \varsigma}
    g_{\bar \sigma, \bar \varsigma}^T \bigg| \varsigma(j_1) \bigg] \Bigg]
    - \EE \Bigg[ \EE \big[ g_{\sigma, \varsigma} \big| \varsigma(j_1) \big] 
    \EE \big[ g_{\bar \sigma, \bar \varsigma}^T \big| \varsigma(j_1) \big]
    \Bigg] = 0.
\end{align*}
Therefore, non-zero terms in Equation (\ref{condUstat:eq:variance_T_4_2}) correspond to the case where there exists $j_3 \neq j_1, j_4 \neq j_1$ such that $\varsigma(j_3) = \bar \varsigma(j_4)$. It is equivalent to
$|\varsigma \cap \bar \varsigma| \geq 2$. We will ignore higher-order terms, i.e. the ones for which $|\varsigma \cap \bar \varsigma| > 2$, as they yield higher powers of $h^{p}$ and are therefore negligible.
Finally, Equation (\ref{condUstat:eq:variance_T_4_2}) becomes
\begin{align*}
    \EE \big[T_{4,2} T_{4,2}^T \big]
    &\simeq \frac{\tilde r_{n,n'}^2}{|\Iknpr|^2 \cdot |\Ikn|^2}
    \sum_{\sigma, \bar \sigma \in \Iknpr}
    \sum_{\substack{\varsigma, \bar \varsigma \in \Ikn \\
    |\varsigma \cap \bar \varsigma| = 2} } \Bigg(
    \EE \Big[ g_{\sigma, \varsigma} g_{\bar \sigma, \bar \varsigma}^T \Big]
    - 2k \EE \bigg[
    \EE \big[ g_{\sigma, \varsigma} \big| \varsigma(i) \big]
    \EE \big[ g_{\bar \sigma, \bar \varsigma}^T \big| \bar \varsigma(\bar i) \big]
    \bigg] \Bigg).
\end{align*}
As before, using change of variables and limited expansions, we can prove that
\begin{align*}
    \frac{\tilde r_{n,n'}^2}{|\Iknpr|^2 \cdot |\Ikn|^2}
    \sum_{\sigma, \bar \sigma \in \Iknpr}
    \sum_{\substack{\varsigma, \bar \varsigma \in \Ikn \\
    |\varsigma \cap \bar \varsigma| = 2} } 
    \EE \Big[ g_{\sigma, \varsigma} g_{\bar \sigma, \bar \varsigma}^T \Big]
    = o(1),
\end{align*}
and similarly for the other term.

\subsection{Convergence of \texorpdfstring{$T_6$}{T6} to \texorpdfstring{$0$}{0}}
\label{section:proof:limit_T_6}

Using Equation (\ref{condUstat:eq:def:term_T_6}), we have $T_6 = T_{6,1} + T_{6,2}$,
where
\begin{align}
    T_{6,1} &:= \frac{\tilde r_{n,n'}}{|\Iknpr|\cdot|\Ikn|}
    \sum_{\sigma \in \Iknpr} \sum_{\varsigma \in \Ikn}
    K_{\sigma, \varsigma}
    \big( g_\varsigma - \theta_\sigma \big)
    \dfrac{N_\sigma - \EE[N_\sigma]}{f_{\Z', \sigma}^2}
    \bar \psibm_\sigma , \\
    T_{6,2} &:= \frac{\tilde r_{n,n'}}{|\Iknpr|\cdot|\Ikn|}
    \sum_{\sigma \in \Iknpr} \sum_{\varsigma \in \Ikn}
    K_{\sigma, \varsigma}
    \big( g_\varsigma - \theta_\sigma \big)
    \dfrac{\EE[N_\sigma] - f_{\Z', \sigma}}{f_{\Z', \sigma}^2}
    \bar \psibm_\sigma.
\end{align}
We first prove that $T_{6,1} = o(1)$. Using Equation (\ref{condUstat:def:normalization_factor_N_k}), we have
\begin{align*}
    T_{6,1} &= \frac{\tilde r_{n,n'}}{|\Iknpr|\cdot|\Ikn|}
    \sum_{\sigma \in \Iknpr} \sum_{\varsigma \in \Ikn} \frac{1}{f_{\Z', \sigma}^2}
    K_{\sigma, \varsigma}
    \big( g_\varsigma - \theta_\sigma \big)
    \big( N_k(\z'_{\sigma(1:k)}) - \EE[ N_k(\z'_{\sigma(1:k)}) ] \big)
    \bar \psibm_\sigma \\
    &= \frac{\tilde r_{n,n'}}{|\Iknpr|\cdot|\Ikn|}
    \sum_{\sigma \in \Iknpr} \sum_{\varsigma \in \Ikn} \frac{1}{f_{\Z', \sigma}^2}
    K_{\sigma, \varsigma}
    \big( g_\varsigma - \theta_\sigma \big)
    \sum_{\nu \in \Ikn} \bigg( 
    \prod_{i=1}^k K_h \big(\Z_{\nu(i)}-\z'_{\sigma(i)} \big)
    - \EE \Big[ \prod_{i=1}^k K_h \big(\Z_{\nu(i)}-\z'_{\sigma(i)} \big) \Big] \bigg)
    \bar \psibm_\sigma \\
    &= \frac{\tilde r_{n,n'}}{|\Iknpr|\cdot|\Ikn|}
    \sum_{\sigma \in \Iknpr} \sum_{\varsigma, \nu \in \Ikn} \frac{1}{f_{\Z', \sigma}^2}
    K_{\sigma, \varsigma}
    \big( g_\varsigma - \theta_\sigma \big)
    \Big( K_{\sigma, \nu} - \EE \big[ K_{\sigma, \nu} \big] \Big)
    \bar \psibm_\sigma.
\end{align*}
The terms for which $|\varsigma \cap \nu| \geq 1$ induce some powers of
$(n h^{p})^{-1}$, and are therefore negligible. We remove them to obtain an equivalent random vector $\bar T_{6,1}$, which is centered.
Therefore it is sufficient to show that its second moment tends to $0$.
\begin{align*}
    \EE \big[ \bar T_{6,1} \bar T_{6,1}^T \big] 
    &= \frac{\tilde r_{n,n'}^2}{|\Iknpr|^2 \cdot |\Ikn|^2}
    \sum_{\sigma, \bar \sigma \in \Iknpr}
    \sum_{\substack{\varsigma, \nu \in \Ikn \\
    \varsigma \cap \nu = \emptyset} }
    \sum_{\substack{\bar \varsigma, \bar \nu \in \Ikn \\
    \bar \varsigma \cap \bar \nu = \emptyset} }
    \frac{\bar \psibm_\sigma}{f_{\Z', \sigma}^2}
    \frac{\bar \psibm_{\bar \sigma}^T}{f_{\Z', \bar \sigma}^2}
    g_{\sigma, \bar \sigma, \varsigma, \bar \varsigma, \nu, \bar \nu}, \\
    g_{\sigma, \bar \sigma, \varsigma, \bar \varsigma, \nu, \bar \nu}
    &:= \EE \bigg[ K_{\sigma, \varsigma}
    \big( g_\varsigma - \theta_\sigma \big)
    \Big( K_{\sigma, \nu} - \EE \big[ K_{\sigma, \nu} \big] \Big)
    K_{\bar \sigma, \bar \varsigma}
    \big( g_{\bar \varsigma} - \theta_{\bar \sigma} \big)
    \Big( K_{\bar \sigma, \bar \nu}
    - \EE \big[ K_{\bar \sigma, \bar \nu} \big] \Big) \bigg].
\end{align*}
The term $g_{\sigma, \bar \sigma, \varsigma, \bar \varsigma, \nu, \bar \nu}$ is $0$ in two cases : if $\nu \cap (\varsigma \cup \bar \varsigma \cup \bar \nu)$ or if $\bar \nu \cap (\varsigma \cup \bar \varsigma \cup \nu)$. This condition can be written as
\begin{align*}
    \emptyset &= \big[ \nu \cap (\bar \varsigma \cup \bar \nu) \big]
    \cup \big[ \bar \nu \cap (\varsigma \cup \nu) \big]
    = (\nu \cup \bar \nu) \cap (\bar \varsigma \cup \bar \nu)
    \cap (\varsigma \cup \nu).
\end{align*}
We deduce that non-zero terms arise only when there exists
$j_1, j_2 \in \{1, \dots, k\}$ such that: $\nu(j_1) = \bar \nu(j_2)$ or $\nu(j_1) = \bar \varsigma(j_2)$ or $\bar \nu(j_1) = \varsigma(j_2)$.
Therefore, we can write $\EE \big[ \bar T_{6,1} \bar T_{6,1}^T \big]
= T_{6,1,1} + T_{6,1,2} + T_{6,1,3}$, where
\begin{align*}
    T_{6,1,1}
    &= \frac{\tilde r_{n,n'}^2}{|\Iknpr|^2 \cdot |\Ikn|^2}
    \sum_{j_1, j_2 = 1}^k
    \sum_{\sigma, \bar \sigma \in \Iknpr}
    \sum_{\substack{\varsigma, \nu \in \Ikn \\
    \varsigma \cap \nu = \emptyset} }
    \sum_{\substack{\bar \varsigma, \bar \nu \in \Ikn \\
    \bar \varsigma \cap \bar \nu = \emptyset, \bar \nu(j_2) = \nu(j_1)} }
    \frac{\bar \psibm_\sigma}{f_{\Z', \sigma}^2}
    \frac{\bar \psibm_{\bar \sigma}^T}{f_{\Z', \bar \sigma}^2}
    g_{\sigma, \bar \sigma, \varsigma, \bar \varsigma, \nu, \bar \nu}, \\
    T_{6,1,2}
    &= \frac{\tilde r_{n,n'}^2}{|\Iknpr|^2 \cdot |\Ikn|^2}
    \sum_{j_1, j_2 = 1}^k
    \sum_{\sigma, \bar \sigma \in \Iknpr}
    \sum_{\substack{\varsigma, \nu \in \Ikn \\
    \varsigma \cap \nu = \emptyset} }
    \sum_{\substack{\bar \varsigma, \bar \nu \in \Ikn \\
    \bar \varsigma \cap \bar \nu = \emptyset,
    \bar \varsigma(j_2) = \nu(j_1)} }
    \frac{\bar \psibm_\sigma}{f_{\Z', \sigma}^2}
    \frac{\bar \psibm_{\bar \sigma}^T}{f_{\Z', \bar \sigma}^2}
    g_{\sigma, \bar \sigma, \varsigma, \bar \varsigma, \nu, \bar \nu}, \\
    T_{6,1,3}
    &= \frac{\tilde r_{n,n'}^2}{|\Iknpr|^2 \cdot |\Ikn|^2}
    \sum_{j_1, j_2 = 1}^k
    \sum_{\sigma, \bar \sigma \in \Iknpr}
    \sum_{\substack{\varsigma, \nu \in \Ikn \\
    \varsigma \cap \nu = \emptyset} }
    \sum_{\substack{\bar \varsigma, \bar \nu \in \Ikn \\
    \bar \varsigma \cap \bar \nu = \emptyset,
    \bar \nu(j_1) = \varsigma(j_2)} }
    \frac{\bar \psibm_\sigma}{f_{\Z', \sigma}^2}
    \frac{\bar \psibm_{\bar \sigma}^T}{f_{\Z', \bar \sigma}^2}
    g_{\sigma, \bar \sigma, \varsigma, \bar \varsigma, \nu, \bar \nu},
\end{align*}
We will prove that $T_{6,1,1} = o(1)$. The two other terms can be treated in a similar way. Because of our assumptions, the terms for which
$\bar \sigma(j_1) \neq \sigma(j_2)$ are zero. This divides the number of possible terms by $n'$.
By using limited expansions as in Equation (\ref{condUstat:eq:expectancy_T_4}), we get that $g_{\sigma, \bar \sigma, \varsigma, \bar \varsigma, \nu, \bar \nu} = O(h^{k\alpha - p})$.
Therefore, we have $T_{6,1,1}
= O \big(\frac{n n' h^{p}}{n n'} h^{k\alpha - p} \big)
= O(h^{k\alpha}) = o(1).$

\mds

Concerning $T_{6,2}$, its variance matrix is given by
\begin{align*}
    Var \big[ T_{6,2} \big] 
    &= \frac{\tilde r_{n,n'}^2}{|\Iknpr|^2 \cdot |\Ikn|^2}
    \sum_{\sigma, \bar \sigma \in \Iknpr}
    \sum_{\varsigma, \bar \varsigma \in \Ikn}
    \dfrac{\EE[N_\sigma] - f_{\Z', \sigma}}{f_{\Z', \sigma}^2}
    \dfrac{\EE[N_{\bar \sigma}] - f_{\Z', \bar \sigma}}{f_{\Z', \bar \sigma}^2}
    \bar \psibm_\sigma \bar \psibm_{\bar \sigma} 
    \bar g_{\sigma, \bar \sigma, \varsigma, \bar \varsigma}, \\
    \bar g_{\sigma, \bar \sigma, \varsigma, \bar \varsigma}
    &:= \EE \bigg[ K_{\sigma, \varsigma} K_{\bar \sigma, \bar \varsigma}
    \big( g_\varsigma - \theta_\sigma \big)
    \big( g_{\bar \varsigma} - \theta_{\bar \sigma} \big) \bigg]
    - \EE \bigg[ K_{\sigma, \varsigma}
    \big( g_\varsigma - \theta_\sigma \big) \bigg]
    \EE \bigg[ K_{\bar \sigma, \bar \varsigma}
    \big( g_{\bar \varsigma} - \theta_{\bar \sigma} \big) \big) \bigg].
\end{align*}
Note that $\bar g_{\sigma, \bar \sigma, \varsigma, \bar \varsigma} = 0$
when $\varsigma \cap \bar \varsigma = \emptyset$.
This divides the number of terms in the sum above by $n$, and imposes that $\sigma \cap \bar \varsigma \neq 0$, which divides the number of terms in the sum above by another $n'$. Finally, limited expansions gives a bound of
$h^{k \alpha - p}$. Summing up all these elements, we obtain
$Var \big[ T_{6,2} \big] = O(\frac{\tilde r_{n,n'}^2}{n n'} h^{k \alpha - p})
= O(h^{k \alpha}) = o(1)$.
Similarly, we get $\EE \big[ T_{6,2} \big] = o(1)$ by a Taylor expansion.

\subsection{Convergence of \texorpdfstring{$T_7$}{T7} to \texorpdfstring{$0$}{0}}
\label{section:proof:limit_T_7}

We recall Equation (\ref{condUstat:eq:def:term_T_7}):
\begin{align*}
    T_7 &= \frac{\tilde r_{n,n'}}{|\Iknpr|\cdot|\Ikn|}
    \sum_{\sigma \in \Iknpr} \sum_{\varsigma \in \Ikn}
    K_{\sigma, \varsigma}
    \big( g_\varsigma - \theta_\sigma \big) T_{7,\sigma}
    \bar \psibm_\sigma, \\
    T_{7, \sigma} &:= \frac{1}{f_{\Z', \sigma}} (1 + \alpha_{7,\sigma})^{-3}
    \left( \dfrac{N_\sigma - f_{\Z', \sigma}}{f_{\Z', \sigma}}
    \right)^2,
    \text{ with } |\alpha_{7,\sigma}| \leq
    \left| \dfrac{N_\sigma - f_{\Z', \sigma}}{f_{\Z', \sigma}} \right|.
\end{align*}
By Lemma~\ref{lemma:bound_Nk_fk} applied to $\z_1 = \z'_{\sigma(1)}, \dots, \z_{n'} = \z'_{\sigma(n')}$, for $\sigma \in \Iknpr$, we get
\begin{align*}
    \PP \bigg( \sup_{\sigma \in \Iknpr} \big| N_\sigma &- f_{\Z', \sigma} \big| \leq \frac{C_{K, \alpha}}{\alpha} h^{\alpha} + t \bigg)
    \geq 1 - 2 \exp \bigg( - \frac{[ n / k] t^2}{h^{-k p} C_1 + h^{-k p} C_2 t} \bigg),
\end{align*}
for any $t > 0$. Therefore,
$\sup_{\sigma \in \Iknpr} |T_{7, \sigma}|
= O_{\PP}(h^{2 \alpha})$ by choosing $t = h^{\alpha / k}$. Then,
\begin{align*}
    |T_7| &\leq \sup_{\sigma \in \Iknpr} |T_{7, \sigma}|
    \frac{\tilde r_{n,n'}}{|\Iknpr|\cdot|\Ikn|}
    \sum_{\sigma \in \Iknpr} \sum_{\varsigma \in \Ikn}
    | K_{\sigma, \varsigma} | \cdot
    \big| g_\varsigma - \theta_\sigma \big| \cdot
    |\bar \psibm_\sigma|.
\end{align*}
The expectation of the double sum is $O(h^\alpha)$, by $\alpha$-order limited expansions. By Markov's inequality, we deduce
\begin{align*}
    T_7 = O_{\PP} \Big( \tilde r_{n,n'}
    \sup_{\sigma \in \Iknpr} |T_{7, \sigma}| h^\alpha) \Big)
    = O_{\PP}(\tilde r_{n,n'} h^{3 \alpha})
    = O_{\PP}\Big( (n n' h^{p + 3 \alpha})^{1/2} \Big),
\end{align*}
therefore $T_7 = o_{\PP}(1)$.

\subsection{Convergence of \texorpdfstring{$T_3$}{T3} to \texorpdfstring{$0$}{0}}
\label{section:proof:limit_T_3}

We have
\begin{align*}
    T_3 := \frac{\tilde r_{n,n'}}{|\Iknpr|}
    \sum_{\sigma \in \Iknpr} \alpha_{3, \sigma} \cdot
    \Big( \hat \theta \big(\z'_{\sigma(1)}, \dots, \z'_{\sigma(k)} \big)
    - \theta \big(\z'_{\sigma(1)}, \dots, \z'_{\sigma(k)} \big) \Big)^2
    \psibm \big(\z'_{\sigma(1)}, \dots, \z'_{\sigma(k)} \big),
\end{align*}
with $\forall \sigma \in \Iknpr$, $|\alpha_{3,\sigma}| \leq C_{\Lambda''}/2$.
Therefore
\begin{align*}
    T_3
    &\lesssim \frac{\tilde r_{n,n'}}{|\Iknpr|}
    \sum_{\sigma \in \Iknpr}
    \Big( \hat \theta \big(\z'_{\sigma(1)}, \dots, \z'_{\sigma(k)} \big)
    - \theta \big(\z'_{\sigma(1)}, \dots, \z'_{\sigma(k)} \big) \Big)^2 \\
    &\lesssim \frac{\tilde r_{n,n'}}{|\Iknpr|} \bigg( \frac{1}{|\Ikn|} \sum_{\varsigma \in \Ikn}
    \dfrac{ K_{\sigma, \varsigma} } { f_{\Z', \sigma} }
    \big( g_\varsigma - \theta_\sigma \big)
    + K_{\sigma, \varsigma} \big( g_\varsigma - \theta_\sigma \big)
    \Big( \dfrac{1} {N_\sigma} - \dfrac{1} {f_{\Z', \sigma}} \Big)
    \bigg)^2 = T_8 + T_9 + T_{10},
\end{align*}
where
\begin{align*}
    T_8
    &:= \frac{\tilde r_{n,n'}}{|\Iknpr| \cdot |\Ikn|^2}
    \sum_{\sigma \in \Iknpr} \sum_{\varsigma, \bar \varsigma \in \Ikn}
    \dfrac{ K_{\sigma, \varsigma} K_{\sigma, \bar \varsigma} }
    { f_{\Z', \sigma}^2 }
    \big( g_\varsigma - \theta_\sigma \big)
    \big( g_{\bar \varsigma} - \theta_\sigma \big), \\
    T_9
    &:= \frac{\tilde r_{n,n'}}{|\Iknpr| \cdot |\Ikn|^2}
    \sum_{\sigma \in \Iknpr} \sum_{\varsigma, \bar \varsigma \in \Ikn}
    \dfrac{ K_{\sigma, \varsigma} K_{\sigma, \bar \varsigma}} { f_{\Z', \sigma} }
    \big( g_\varsigma - \theta_\sigma \big)
    \big( g_{\bar \varsigma} - \theta_\sigma \big)
    \Big( \dfrac{1} {N_\sigma} - \dfrac{1} {f_{\Z', \sigma}} \Big), \\
    T_{10}
    &:= \frac{\tilde r_{n,n'}}{|\Iknpr| \cdot |\Ikn|^2}
    \sum_{\sigma \in \Iknpr} \sum_{\varsigma, \bar \varsigma \in \Ikn}
    K_{\sigma, \varsigma} K_{\sigma, \bar \varsigma}
    \big( g_\varsigma - \theta_\sigma \big)
    \big( g_{\bar \varsigma} - \theta_\sigma \big)
    \Big( \dfrac{1} {N_\sigma} - \dfrac{1} {f_{\Z', \sigma}} \Big)^2.
\end{align*}

We show that $T_8 = o(1)$. The two other terms can be treated in a similar way.
\begin{align*}
    \EE \big[|T_8| \big]
    &= \EE \bigg[ \frac{\tilde r_{n,n'}}{|\Iknpr| \cdot |\Ikn|^2}
    \sum_{\sigma \in \Iknpr} \sum_{\varsigma, \bar \varsigma \in \Ikn}
    \dfrac{ |K_{\sigma, \varsigma} K_{\sigma, \bar \varsigma}| }
    { f_{\Z', \sigma}^2 }
    \big| g_\varsigma - \theta_\sigma \big| \cdot
    \big| g_{\bar \varsigma} - \theta_\sigma \big| \bigg] \\
    &= \frac{\tilde r_{n,n'}}{|\Iknpr| \cdot |\Ikn|^2}
    \sum_{\sigma \in \Iknpr} \sum_{\varsigma, \bar \varsigma \in \Ikn}
    \int \frac{\prod_{i=1}^k  \Big| K_h \big(\z_{\varsigma(i)}-\z'_{\sigma(i)} \big)
    K_h \big(\z_{\bar \varsigma(i)}-\z'_{\sigma(i)} \big) \Big|} {f_{\Z', \sigma}^2} \\
    & \hspace{1cm} \cdot 
    \Big| g \big(\x_{\varsigma(1:k)} \big) - \theta_\sigma \Big| 
    \Big| g \big(\x_{\bar \varsigma(1:k)} \big) - \theta_\sigma \Big| 
    \prod_{i \, \in \, \varsigma(1:k) \, \cup \, \bar \varsigma(1:k)}
    f_{\X, \Z} (\x_i, \z_i) d\mu(\x_i) d\z_i.
\end{align*}
Note that terms for which $\varsigma \neq \bar \varsigma \in \Iknpr$ are zero, because the $\z'_i$ are distinct and because of our Assumption~\ref{condUstat:assumpt:asymptNorm_joint}\ref{condUstat:assumpt:kernel_separation_Zprime}.
Therefore, we get
\begin{align*}
    \EE \big[|T_8| \big]
    &= \frac{\tilde r_{n,n'}}{|\Iknpr| \cdot |\Ikn|^2}
    \sum_{\sigma \in \Iknpr} \sum_{\varsigma \in \Ikn}
    \int \frac{\prod_{i=1}^k K_h \big(\z_{\varsigma(i)}-\z'_{\sigma(i)} \big)^2} {f_{\Z', \sigma}^2} 
    \Big( g \big(\x_{\varsigma(1:k)} \big) - \theta_\sigma \Big)^2
    \prod_{i \, \in \, \varsigma(1:k)} f_{\X, \Z} (\x_i, \z_i) d\mu(\x_i) d\z_i \\
    &= \frac{\tilde r_{n,n'}}{|\Iknpr| \cdot |\Ikn|} \sum_{\sigma \in \Iknpr}
    \int \frac{\prod_{i=1}^k K_h \big(\z_i-\z'_{\sigma(i)} \big)^2} {f_{\Z', \sigma}^2} 
    \Big( g \big(\x_{\varsigma(1:k)} \big) - \theta_\sigma \Big)^2
    \prod_{i=1}^k f_{\X, \Z} (\x_i, \z_i) d\mu(\x_i) d\z_i \\
    &= \frac{\tilde r_{n,n'} h^{- k p} }{|\Iknpr| \cdot |\Ikn|}
    \sum_{\sigma \in \Iknpr}
    \int \frac{\prod_{i=1}^k K \big( \t_i \big)^2} {f_{\Z', \sigma}^2} 
    \Big( g \big(\x_{\varsigma(1:k)} \big) - \theta_\sigma \Big)^2
    \prod_{i=1}^k f_{\X, \Z} (\x_i, \z'_{\sigma(i)} + h \t_i) d\mu(\x_i) d\z_i \\
    &= O \bigg( \frac{\tilde r_{n,n'} h^{- k p} }{|\Ikn|} \bigg)
    = O \bigg( \Big( \frac{n \times n' \times h^{(1 - k) p} }{|\Ikn|^2} \Big)^{1/2} \bigg) = o(1). \; \Box
\end{align*}

\end{document}